\documentclass{amsart}%
\usepackage{amsmath}
\usepackage{amsfonts}
\usepackage{amssymb}
\usepackage{graphicx}%
\providecommand{\U}[1]{\protect\rule{.1in}{.1in}}
\newtheorem{theorem}{Theorem}
\theoremstyle{plain}

\newtheorem{condition}{Condition}

\newtheorem{corollary}{Corollary}

\newtheorem{lemma}{Lemma}

\newtheorem{proposition}{Proposition}
\newtheorem{remark}{Remark}

\numberwithin{equation}{section}
\begin{document}
\title[$q-$Wiener process]{$q-$Wiener and $(\alpha,q)-$ Ornstein--Uhlenbeck processes. A generalization
of known processes}
\author{Pawe\l \ J. Szab\l owski}
\address{Department of Mathematics and Information Sciences\\
Warsaw University of Technology\\
pl. Politechniki 1\\
00-661 Warszawa, Poland}
\email{pszablowski@elka.pw.edu.pl; pawel.szablowski@gmail.com}
\thanks{The author is grateful to an unknown referee for pointing out subtle questions
concerning some proofs, hasty thoughtlessness in posing open questions and
numerous misprints.}
\date{February, 2009}

\begin{abstract}
We collect, scattered through literature, as well as we prove some new
properties of two Markov processes that in many ways resemble Wiener and
Ornstein--Uhlenbeck processes. Although processes considered in this paper
were defined either in non-commutative probability context or through
quadratic harnesses we define them once more as so to say 'continuous time '
generalization of a simple, symmetric, discrete time process satisfying simple
conditions imposed on the form of its first two conditional moments. The
finite dimensional distributions of the first one (say $\mathbf{X=}\left(
X_{t}\right)  _{t\geq0}$ called $q-$Wiener) depends on one parameter
$q\in(-1,1]$ and of the second one (say $\mathbf{Y=}\left(  Y_{t}\right)
_{t\in\mathbb{R}}$ called $(\alpha,q)-$ Ornstein--Uhlenbeck) on two parameters
$\left(  \alpha,q\right)  \in(0,\infty)\times(-1,1]$. The first one resembles
Wiener process in the sense that for $q=1$ it is Wiener process but also that
for $\left\vert q\right\vert <1$ and $\forall n\geq1:$ $t^{n/2}H_{n}\left(
X_{t}/\sqrt{t}|q\right)  ,$ where $\left(  H_{n}\right)  _{n\geq0}$ are the so
called $q-$Hermite polynomials, are martingales. It does not have however
neither independent increments not allows continuous sample path modification.
The second one resembles Ornstein--Uhlenbeck process. For $q=1$ it is a
classical OU process. For $\left\vert q\right\vert <1$ it is also stationary
with correlation function equal to $\exp(-\alpha|t-s|)$ and has many
properties resembling those of its classical version. We think that these
process are fascinating objects to study posing many interesting, open questions.

\end{abstract}
\maketitle

\section{Introduction}

As announced in the abstract, we are going to define two time-continuous
families of Markov processes. One of them will resemble Wiener process and the
other Ornstein--Uhlenbeck (OU) process. They will be indexed (apart from time
parameter) by additional parameter $q\in(-1,1].$ In the case $q\allowbreak
=\allowbreak1$ these processes are classical Wiener and OU process. Of course
for OU process there will be additional parameter $\alpha\allowbreak
>\allowbreak0$ responsible for covariance function of the process. For
$q\in(-1,1)$ both processes will assume values in a compact space :
$(\alpha,q)-$OU process on $[-\frac{2}{\sqrt{1-q}},\frac{2}{\sqrt{1-q}}]$ ,
while for $q-$Wiener process $\left(  X_{t}\right)  _{t\geq0}$ we will have:
$X_{t}\in\lbrack-2\sqrt{\frac{t}{1-q}},2\sqrt{\frac{t}{1-q}}]$. One
dimensional probabilities and transitional probabilities of these processes
will be given explicitly. Moreover two families of polynomials, orthogonal
with respect to these measures, will also be presented. Some properties of
conditional expectations given the past and also past and the future will be
exposed. Martingale properties, as well as some properties of sample path of
these processes will be described. Thus quite detailed knowledge concerning
these processes will be presented.

The processes that we are going to reintroduce have appeared already in $1997$
in an excellent paper \cite{Bo} as an offspring and a particular case of some
non-commutative probability model. Some of the properties of these processes
particularly those associated with martingale behavior of some functions of
these processes were also discussed in this paper. Since $1997$ there appeared
couple of papers on the properties of $q-$Gaussian distributions. See e.g.
\cite{BryBo}, \cite{Bry3}, \cite{Ansh}.

There is also a different path of research followed by W\l odek Bryc, Wojtek
Matysiak and Jacek Weso\l owski see e.g. \cite{bryc05}, \cite{BryWe},
\cite{BryMaWe}. Their starting point is a process with continuous time that
satisfies several (exactly $5$) conditions on covariance function and on the
first and the second conditional moments. Those are the so called quadratic
harnesses characterized by $5$ parameters. Under resulting $5$ assumptions
they proved that these processes are Markov and also stated several properties
of the families of polynomials that orthogonalize transitional and one
dimensional probabilities. They gave several examples illustrating developed
theory. One of the processes considered by them is the so called $q-$Brownian
process. $4$ of $5$ possible parameters are equal to zero and the fifth one
can be identified with parameter $q$ considered in this paper. As far as the
one dimensional probabilities and the transitional probabilities are concerned
this process is identical with $q-$Wiener process introduced and analyzed in
this paper. They did not however work on the properties of the $q-$Brownian
process. It appeared as a by-product of their interest in quadratic harnesses.
Bryc Matysiak and Weso\l owski were mostly interested in the general problem
of existence of quadratic harnesses. That is why $(\alpha,q)-$OU process have
not appeared in their works.

What we are aiming to do is to reintroduce these processes via certain
discrete time one dimensional, time symmetric random process (1TSP) by the so
to say "continuation of time" or may be more precisely as processes that
sampled at certain discrete moments have the properties of 1TSP. This discrete
time process was defined in a purely classical probability context. Moreover
it is very simple and intuitive. Its simplicity surprised Richard Askey "that
it has no $q$ in the statement of the problem" as he puts it in the forward to
\cite{IA}. Thus these 1TSP can be simple models of some phenomenons observed
in the recently intensively developing $q$-series theory. 1TSP first appeared
in 2000 in \cite{bryc1} and have been studied in detail recently see e.g.
\cite{bryc2}, \cite{bms}, \cite{matszab}, \cite{matszab2}.

Besides we derive these processes under fewer assumptions (than one would need
while following quadratic harnesses path) on the first and the second
conditional moments (we need only $2$). The construction is also different. As
mentioned earlier our starting point is a discrete time 1TSP. $q-$Wiener
process is obtained as continuous time transformation of the process that we
call $(\alpha,q)-$OU process (like in the classical i.e. $q\allowbreak
=\allowbreak1$ case) which, on its side, is obtained as continuous time
generalization of the discrete time process (1TSP). Besides we list many
properties of $q-$Wiener process that justify its name (are sort of $q-$
analogies of well known martingale properties of Wiener process). Some of
these properties can be derived from the definition of $q-$Brownian process
presented in \cite{bryc05}. They are not however stated explicitly.

That is why the first section will be dedicated to definition of 1TSP's and
recollection of their basic properties. Since our approach is totally
commutative $q-$ by no means is an operator. It is a number parameter. Yet we
are touching the $q-$ series theory and the special functions.

We think that processes presented in this paper are the fascinating objects to
study. As mentioned before, for $\left\vert q\right\vert <1:$ $q-$Wiener
process has many properties similar to ordinary Brownian motion, but is not an
independent increment process. Besides we present a few properties of
trajectories of the processes discussed in this paper. Properties that,
although can be relatively easily deduced, were never stated in the above
mentioned papers where these processes appeared first.

The paper is organized as follows. As we mentioned, in the second section we
recall definition and summarize basic properties of 1TSP's. The next section
is still dedicated mainly to certain auxiliary properties of the discrete time
1TSP's that are necessary to perform our construction. In the fourth section
we introduce continuous version of 1TSP ($(\alpha,q)-$OU process) and prove
its existence. Then we define $q-$Wiener as a continuous time transformation
of $(\alpha,q)-$OU process. Later parts of this section are devoted to
presentation of the two processes and listing or proving some of their
properties. We indicate their connection with an emerging quadratic harness theory.

We point out here where is the mistake causing that Weso\l owski's martingale
characterization of the Wiener process contained in \cite{Wes} is not true.
The fact that is not true was already known to Weso\l owski (see e.g.
\cite{bryc05} where processes denying his characterization are pointed out).

The fifth section presents some obvious open problems that come to mind almost
directly. The last (sixth) section contains lengthy proofs of the results from
the previous sections.

It is known that 1TSP exist with parameter $q>1$ (see \cite{Szab}). It's
transition distribution is then discrete. We show in this paper that
$q-$Wiener and $(\alpha,q)-$Ornstein--Uhlenbeck processes do not exist for
$q>1.$

\section{One dimensional time symmetric random processes\label{wlasnosci}}

By 1TSP's we mean square integrable random field $\mathbf{X=}\left\{
X_{n}\right\}  _{n\in\mathbb{Z}}$ indexed by the integers, with non-singular
all covariance matrices and constant first two moments, that satisfy the
following two sets of conditions :%

\begin{equation}
\exists a,b\in\mathbb{R};\forall n\in\mathbb{Z}:\mathbb{E}\left(
X_{n}|\mathcal{F}_{\neq n}\right)  =a\left(  X_{n-1}+X_{n+1}\right)  +b,~a.s.
\label{bryc1}%
\end{equation}
and
\begin{subequations}
\label{bryc2}%
\begin{align}
\exists A,B,C  &  \in\mathbb{R};\forall n\in\mathbb{Z}:\mathbb{E}\left(
X_{n}^{2}|\mathcal{F}_{\neq n}\right)  =\label{bryc2_1}\\
&  A\left(  X_{n-1}^{2}+X_{n+1}^{2}\right)  +BX_{n-1}X_{n+1}+D\left(
X_{n-1}+X_{n+1}\right)  +C,~a.s., \label{bryc2_2}%
\end{align}
where $\mathcal{F}_{\neq m}:=\sigma\left(  X_{k}:k\neq m\right)  .$

Let us define also $\sigma-$algebras $\mathcal{F}_{\leq m}:=\sigma\left(
X_{k}:k\leq m\right)  ,$ $\mathcal{F}_{\geq m}:=\sigma\left(  X_{k}:k\geq
m\right)  ,$and $\mathcal{F}_{\leq m,\geq j}:=\sigma\left(  X_{k}:k\leq m\vee
k\geq j\right)  $.

Non-singularity of covariance matrices implies that all random variables
$X_{n}$ are non-degenerate and there is no loss of generality in assuming that
$\mathbb{E}X_{k}=0$ and $\mathbb{E}X_{k}^{2}=1$ for all $k\in\mathbb{Z},$
which implies $b=0.$

It has been shown in \cite{matszab} that (\ref{bryc1}) implies $L_{2}%
$-stationarity (stationarity in the wider sense) of $\mathbf{X.}$ Since the
case $\rho:=\operatorname*{corr}\left(  X_{0},X_{1}\right)  =0$ contains
sequences of independent random variables (which satisfy (\ref{bryc1}) and
(\ref{bryc2}) but can have arbitrary distributions), we shall exclude it from
the considerations. Observe that non-singularity of the covariance matrices
implies $\left\vert \rho\right\vert <1.$ By Theorem 3.1 from \cite{bryc1} (see
also Theorems 1 and 2 in \cite{matszab}), we have $\operatorname*{corr}\left(
X_{0},X_{k}\right)  =\rho^{\left\vert k\right\vert }.$ Moreover the one-sided
regressions are linear%
\end{subequations}
\begin{equation}
\mathbb{E}\left(  X_{m}|\mathcal{F}_{\leq0}\right)  =\rho^{m}X_{0}%
=\mathbb{E}\left(  X_{-m}|\mathcal{F}_{\geq0}\right)  ~;~m\geq1.
\label{onesidedregr}%
\end{equation}
It turns out that parameters $a,\rho$, $A,$ $B,$ $C$ are related to one
another. In \cite{bryc1} and \cite{matszab2} it was shown that one can
redefine parameters by introducing new parameter $q\allowbreak=\allowbreak
\frac{\rho^{4}+B\left(  \rho+\frac{1}{\rho}\right)  ^{2}-1}{1+\rho^{4}\left(
B\left(  \rho+\frac{1}{\rho}\right)  ^{2}-1\right)  }$ and express parameters
$A,$ $B,$ $C$ with the help of $\rho$ and $q$ only in the following manner:
\begin{equation}
A=\frac{\rho^{2}\left(  1-q\rho^{2}\right)  }{\left(  \rho^{2}+1\right)
\left(  1-q\rho^{4}\right)  },B=\frac{\rho^{2}\left(  1-\rho^{2}\right)
\left(  1+q\right)  }{\left(  \rho^{2}+1\right)  \left(  1-q\rho^{4}\right)
},C=\frac{\left(  1-\rho^{2}\right)  ^{2}}{1-q\rho^{4}}. \label{_B}%
\end{equation}

We can rephrase and summarize the results of \cite{bryc1} and \cite{matszab}
in the following way. Each 1TSP is characterized by two parameters $\rho$ and
$q.$ For $q$ outside the set $[-1,1]\cup\left\{  1/\rho^{2/n}:n\in
\mathbb{N}\right\}  $ 1TSP's do not exist. For $0<\left\vert \rho\right\vert
<1$ and $q\in\lbrack-1,1]$ 1TSP's exist and all their finite dimensional
distribution are uniquely determined and known. Also it follows from
\cite{bryc1} and \cite{bryc2} that 1TSP are Markov processes. The case
$0<\left\vert \rho\right\vert <1$ and $q\in\left\{  1/\rho^{2/n}%
:n\in\mathbb{N}\right\}  $ is treated in \cite{Szab}.

1TSP with $0<\left\vert \rho\right\vert <1$ and $q\in(-1,1]$ we will call
\emph{regular }.

We adopt notation traditionally used in '$q-$series theory': $\left(
a;q\right)  _{0}\allowbreak=\allowbreak1,$
\[
\left(  a;q\right)  _{n}\allowbreak=\allowbreak%
{\displaystyle\prod\limits_{i=0}^{n-1}}
(1-aq^{i}),\left(  a_{1},\ldots,a_{k};q\right)  _{n}\allowbreak=\allowbreak
\prod_{i=1}^{k}\left(  a_{i};q\right)  _{n},
\]
$\left[  0\right]  _{q}\allowbreak=\allowbreak0,$ $\left[  n\right]
_{q}\allowbreak=\allowbreak1+\ldots q^{n-1},$ $n\geq1,$ $\left[  0\right]
_{q}!\allowbreak=\allowbreak1,$ $\left[  n\right]  _{q}!\allowbreak
=\allowbreak%
{\displaystyle\prod\limits_{i=1}^{n}}
\left[  i\right]  _{q},$
\[%
\genfrac{[}{]}{0pt}{}{n}{k}%
_{q}\allowbreak=\allowbreak\left\{
\begin{array}
[c]{ccc}%
\frac{\left[  n\right]  _{q}!}{\left[  k\right]  _{q}!\left[  n-k\right]
_{q}} & when & 0\leq k\leq n\\
0 & when & k>n
\end{array}
\right.  .
\]
Notice that we have: $\left(  q;q\right)  _{n}\allowbreak=\allowbreak\left(
1-q\right)  ^{n}\left[  n\right]  _{n}!,$ $%
\genfrac{[}{]}{0pt}{}{n}{k}%
_{q}\allowbreak=\allowbreak\frac{\left(  q;q\right)  _{n}}{\left(  q;q\right)
_{k}\left(  q;q\right)  _{n-k}}$. As it is customary in $q-$series theory we
will often abbreviate $\left(  a;q\right)  _{n}$and $\left(  a_{1}%
,\ldots,a_{k};q\right)  _{n}$ to $\left(  a\right)  _{n}$ and $\left(
a_{1},\ldots,a_{k}\right)  _{n}$ if it will not cause misunderstanding.

To complete recollection of basic properties of 1TSP let us introduce the so
called $q-$Hermite polynomials $\left\{  H_{n}\left(  x|q\right)  \right\}
_{n\geq-1}$ defined by the following recurrence:
\begin{equation}
\forall n\geq0:xH_{n}\left(  x|q\right)  =H_{n+1}(x|q)+\left[  n\right]
_{q}H_{n-1}\left(  x|q\right)  , \label{Her}%
\end{equation}
with $H_{-1}\left(  x|q\right)  =0,~H_{0}\left(  x|q\right)  =1.$

\begin{remark}
Comparing initial values and $3$-term recurrences (see e.g. \cite{Andrews1999}%
) one can easily notice that $\left\{  H_{n}\left(  x|1\right)  \right\}
_{n\geq-1}$ are the so called probabilistic Hermite polynomials (i.e.
orthogonal with respect to the measure with density $\exp(-x^{2}/2)/\sqrt
{2\pi})$ ), while for $\forall n\geq-1$ $H_{n}\left(  x|0\right)
\allowbreak=\allowbreak U_{n}\left(  x/2\right)  $ where $\left\{
U_{n}\right\}  $ are the so called Chebyshev polynomials of the second kind
i.e. polynomials orthogonal with respect to the measure with density
$\sqrt{1-x^{2}}/\pi.$
\end{remark}

\begin{enumerate}
\item W. Bryc in \cite{bryc1} has shown that there exist stationary
distribution of $\mathbf{X}$ and that:
\[
\forall n\in\mathbb{Z},k,i\geq1:\mathbb{E}\left(  H_{k}\left(  X_{n}|q\right)
|\mathcal{F}_{\leq n-i}\right)  =\rho^{ki}H_{k}\left(  X_{n-i}|q\right)
,~a.s.
\]

\end{enumerate}

and also that $\left\{  H_{n}\left(  x|q\right)  \right\}  _{n\geq-1}$ are
orthogonal polynomials of the stationary distribution.

The case $q\allowbreak=\allowbreak-1$ is equivalent to $B=0$ and leads to
marginal symmetric distribution concentrated on $\left\{  -1,1\right\}  $.

We will concentrate thus on the case $q\in(-1,1]$.

Let $I_{A}\left(  x\right)  \allowbreak=\allowbreak\left\{
\begin{array}
[c]{ccc}%
1 & if & x\in A\\
0 & if & x\notin A
\end{array}
\right.  $ denote index function of the set $A$ and define also
\begin{equation}
S\left(  q\right)  \allowbreak=[-2/\sqrt{1-q},2/\sqrt{1-q}]:\text{for }%
q\in(-1,1)\text{ and }S\left(  1\right)  \allowbreak=\allowbreak\mathbb{R}.
\label{S_q}%
\end{equation}

It turns out that the stationary distribution of 1TSP has for $q\in(-1,1)$ the
density given by
\begin{equation}
f_{N}(x|q)=\frac{\sqrt{1-q}\left(  q\right)  _{\infty}}{2\pi\sqrt
{4-(1-q)x^{2}}}\prod_{k=0}^{\infty}\left(  (1+q^{k})^{2}-(1-q)x^{2}%
q^{k}\right)  I_{S\left(  q\right)  }\left(  x\right)  , \label{q_norm}%
\end{equation}
while for $q\allowbreak=\allowbreak1$ it is equal to
\begin{equation}
f_{N}\left(  x|1\right)  =\frac{1}{\sqrt{2\pi}}\exp\left(  -\frac{x^{2}}%
{2}\right)  ,~x\in\mathbb{R}, \label{Norm1}%
\end{equation}
In particular for $q\allowbreak=\allowbreak0$ we have%
\[
f_{N}\left(  x|0\right)  \allowbreak=\allowbreak\frac{\sqrt{4-x^{2}}}{2\pi}.
\]
W. Bryc in \cite{bryc1} has also found the density of the conditional
distribution $X_{n}|X_{n-1}=y$ and later W. Bryc, W. Matysiak and P.
J.\ Szab\l owski in \cite{bms} have found orthogonal polynomials of this
conditional distribution. Namely it turned out that this distribution has for
$q\in(-1,1)$ and $\allowbreak y\allowbreak\in\allowbreak S\left(  q\right)  $
density of the form
\begin{gather}
f_{CN}\left(  x|y,\rho,q\right)  =\frac{\sqrt{1-q}\left(  \rho^{2},q\right)
_{\infty}}{2\pi\sqrt{4-(1-q)x^{2}}}\times\label{w_kowa}\\
\prod_{k=0}^{\infty}\frac{\left(  (1+q^{k})^{2}-(1-q)x^{2}q^{k}\right)
}{(1-\rho^{2}q^{2k})^{2}-(1-q)\rho q^{k}(1+\rho^{2}q^{2k})xy+(1-q)\rho
^{2}(x^{2}+y^{2})q^{2k}}I_{S\left(  q\right)  }\left(  x\right)  ,\nonumber
\end{gather}
and that polynomials $\left\{  P_{n}\left(  x|y,\rho,q\right)  \right\}
_{n\geq-1}$ defined by
\begin{equation}
\forall n\geq0:P_{n+1}(x|y,q,\rho)=(x-\rho yq^{n})P_{n}(x|y,q,\rho
)-(1-\rho^{2}q^{n-1})\left[  n\right]  _{q}P_{n-1}(x|y\,q,\rho),
\label{alsalamy}%
\end{equation}
with $P_{-1}\left(  x|y,\rho,q\right)  =0,~P_{0}\left(  x|y,q\,,\rho\right)
=1$, are orthogonal with respect to the measure defined by the density
(\ref{w_kowa}). We will call polynomials $\left\{  P_{n}\right\}  $
Al-Salam--Chihara (briefly ASC).

To support intuition let us notice that both densities $f_{N}$ and $f_{CN}$
are bounded. More precisely we have the following easy remark giving bounds
for the both considered densities. Let us remark also that in
\cite{Szablowski2009} are presented bounds for $f_{N}$ more subtle than the
ones given below.

\begin{remark}
i) For $\left\vert q\right\vert <1,~\forall x\in S\left(  q\right)
:f_{N}\left(  x|q\right)  \allowbreak\leq\allowbreak\frac{\sqrt{1-q}}{\pi
}\left(  q\right)  _{\infty}\left(  -\left\vert q\right\vert \right)
_{\infty}^{2},$ and $f_{N}\left(  \frac{\pm2}{\sqrt{1-q}}|q\right)
\allowbreak=\allowbreak0,$

ii) For $\left\vert q\right\vert <1,~\forall x,y\in S\left(  q\right)
:0<C\left(  y,\rho,q\right)  \leq\frac{f_{CN}\left(  x|y,\rho,q\right)
}{f_{N}\left(  x|q\right)  }\leq\frac{\left(  \rho^{2}\right)  _{\infty}%
}{\left(  \rho\right)  _{\infty}^{4}}.$
\end{remark}

\begin{proof}
i) Follows the fact that for $x\in S\left(  q\right)  :\left(  (1+q^{k}%
)^{2}-(1-q)x^{2}q^{k}\right)  \allowbreak\leq\allowbreak\left(  1+\left\vert
q\right\vert ^{k}\right)  $ and $\sqrt{4-x^{2}}\leq2.$ ii) was proved in
\cite{szab2011} Proposition 1,vii.
\end{proof}

Again we have two special simple cases presented in the following remark.

\begin{remark}
$\forall n\geq-1$ : $P_{n}\left(  x|y,\rho,1\right)  \allowbreak=\allowbreak
H_{n}\left(  \frac{x-\rho y}{\sqrt{1-\rho^{2}}}\right)  \left(  1-\rho
^{2}\right)  ^{n/2},$ $P_{n}\left(  x|y,\rho,0\right)  \allowbreak=\allowbreak
U_{n}\left(  x/2\right)  -\rho yU_{n}\left(  x/2\right)  \allowbreak
+\allowbreak\rho^{2}U_{n}\left(  x/2\right)  .$
\end{remark}

\begin{proof}
There are many proofs of these simple facts. One of the simplest can be found
in say \cite{Szablowski2010(1)}.
\end{proof}

We will need several properties of polynomials $\left\{  H_{n}\right\}  $ and
$\left\{  P_{n}\right\}  .$ Most of these properties can be found in \cite{IA}
and some in \cite{bms}. We will collect these properties in the following Lemma.

\begin{lemma}
\label{Hip}Assume $\left\vert q\right\vert <1.$ i) For $n,m\geq0:$%
\[
\int_{S\left(  q\right)  }H_{n}\left(  x|q\right)  H_{m}\left(  x|q\right)
f_{N}\left(  x|q\right)  dx\allowbreak=\allowbreak\left\{
\begin{array}
[c]{ccc}%
0 & when & n\neq m\\
\left[  n\right]  _{q}! & when & n=m
\end{array}
\right.  .
\]

ii) For $n\geq0:$%
\[
\int_{S\left(  q\right)  }H_{n}\left(  x|q\right)  f_{CN}\left(
x|y,\rho,q\right)  dx=\rho^{n}H_{n}\left(  y|q\right)  ,
\]

iii) For $n,m\geq0:$%
\[
\int_{S\left(  q\right)  }P_{n}\left(  x|y,\rho,q\right)  P_{m}\left(
x|y,\rho,q\right)  f_{CN}\left(  x|y,\rho,q\right)  dx\allowbreak
\newline=\allowbreak\left\{
\begin{array}
[c]{ccc}%
0 & when & n\neq m\\
\left(  \rho^{2}\right)  _{n}\left[  n\right]  _{q}! & when & n=m
\end{array}
.\right.
\]

iv)
\[
\int_{S\left(  q\right)  }f_{CN}\left(  x|y,\rho_{1},q\right)  f_{CN}\left(
y|z,\rho_{2},q\right)  dy=f_{CN}\left(  x|z,\rho_{1}\rho_{2},q\right)  .
\]

v) $\max_{x\in S\left(  q\right)  }\left\vert H_{n}\left(  x|q\right)
\right\vert \leq W_{n}\left(  q\right)  /\left(  1-q\right)  ^{n/2}$, where
$W_{n}\left(  q\right)  \allowbreak=\allowbreak\sum_{i=0}^{n}%
\genfrac{[}{]}{0pt}{}{n}{i}%
_{q}.$

vi) $\sum_{i=0}^{\infty}\frac{W_{i}\left(  q\right)  t^{i}}{\left(  q\right)
_{i}}\allowbreak=\allowbreak\frac{1}{\left(  t\right)  _{\infty}^{2}}$ and
$\sum_{i=0}^{\infty}\frac{W_{i}^{2}\left(  q\right)  t^{i}}{\left(  q\right)
_{i}}\allowbreak=\allowbreak\frac{\left(  t^{2}\right)  _{\infty}}{\left(
t\right)  _{\infty}^{3}}$ absolutely for $|t|,\left\vert q\right\vert
<1,$\newline where $W_{i}\left(  q\right)  $ is defined in v).

vii) For $(1-q)x^{2}\leq4$ and $\forall(1-q)t^{2}<1\allowbreak:\allowbreak$%
\[
\varphi\left(  x,t|q\right)  \allowbreak\overset{df}{=}\allowbreak\sum
_{i=0}^{\infty}\frac{t^{i}}{\left[  i\right]  _{q}!}H_{n}\left(  x|q\right)
\allowbreak=\prod_{k=0}^{\infty}\left(  1-\left(  1-q\right)  xtq^{k}+\left(
1-q\right)  t^{2}q^{2k}\right)  ^{-1},
\]
convergence is absolute and uniform in $x$. Moreover $\varphi\left(
x,t|q\right)  $ is nonnegative and \newline$\int_{S\left(  q\right)  }%
\varphi\left(  x,t|q\right)  f_{N}\left(  x|q\right)  dx\allowbreak
=\allowbreak1$. $\allowbreak$

viii) For $(1-q)\max(x^{2},y^{2})\leq4,$ $\left\vert \rho\right\vert <1$ and
$\forall(1-q)t^{2}<1\allowbreak:\allowbreak$%
\[
\tau\left(  x,t|y,\rho,q\right)  \allowbreak=\sum_{i=0}^{\infty}\frac{t^{i}%
}{\left[  i\right]  _{q}!}P_{n}\left(  x|y,\rho,q\right)  =\allowbreak
\prod_{k=0}^{\infty}\frac{\left(  1-\left(  1-q\right)  \rho ytq^{k}+\left(
1-q\right)  \rho^{2}t^{2}q^{2k}\right)  }{\left(  1-\left(  1-q\right)
xtq^{k}+\left(  1-q\right)  t^{2}q^{2k}\right)  },
\]
convergence is absolute and uniform in $x$. Moreover $\tau\left(
x,t|\theta,\rho,q\right)  $ is nonnegative and $\int_{S\left(  q\right)  }%
\tau\left(  x,t|y,\rho,q\right)  f_{CN}\left(  x|y,\rho,q\right)
dx\allowbreak=\allowbreak1$.

ix) For $(1-q)\max(x^{2},y^{2})\leq2,$ $\left\vert \rho\right\vert <1$%
\[
f_{CN}\left(  x|y,\rho,q\right)  \allowbreak\allowbreak=\allowbreak
f_{N}\left(  x|q\right)  \sum_{n=0}^{\infty}\frac{\rho^{n}}{[n]_{q}!}%
H_{n}(x|q)H_{n}(y|q)
\]
and convergence is absolute and uniform in $x$ and $y.$
\end{lemma}

\begin{proof}
Since the more popular are slightly modified polynomials $H_{n}$ namely
polynomials $h_{n}\left(  x|q\right)  \allowbreak=\allowbreak H_{n}\left(
\frac{2}{\sqrt{1-q}}|q\right)  /\left(  1-q\right)  ^{n/2};$ $x\in
\lbrack-1,1]$ (called continuous $q-$Hermite polynomials) and $p_{n}\left(
x|y,\rho,q\right)  \allowbreak=\allowbreak P_{n}\left(  \frac{2x}{\sqrt{1-q}%
}|\frac{2y}{\sqrt{1-q}},\rho,q\right)  $, many properties of $q-$Hermite and
Al-Salam--Chihara polynomials are formulated in terms of $h_{n}$ and $p_{n}.$
i) It is formula 13.1.11 of \cite{IA} with an obvious modification for the
polynomials $H_{n}$ instead of $h_{n}$ and normalized weight function (i.e.
$f_{N})$ ii) Exercise 15.7 of \cite{IA} also in \cite{bryc1}, iii) Formula
15.1.5 of \cite{IA} with obvious modification for polynomials $p_{n}$ instead
of $P_{n}$ and normalized weight function (i.e. $f_{CN}),$ iv) see (2.6) of
\cite{bms}. v) and vi) Exercise 12.2(b) and 12.2(c) of \cite{IA}. vii)-viii)
follow v) and vi). Besides non-negativity of $\varphi$ and $\tau$ are trivial
and follow formulae $1-\left(  1-q\right)  xtq^{k}+\left(  1-q\right)
t^{2}q^{2k}\allowbreak=\allowbreak(1-q)(tq^{k}-x/2)^{2}\allowbreak
+\allowbreak1-(1-q)x^{2}/4$ and $1-\left(  1-q\right)  \rho ytq^{k}+\left(
1-q\right)  \rho^{2}t^{2}q^{2k}\allowbreak=\allowbreak(1-q)\rho^{2}%
(q^{k}t-y/(2\rho))^{2}+1-(1-q)y^{2}/4$. Values of integrals follow i) and
iii). ix) is the famous Poisson-Mehler expansion formula. It has many proofs
presented e.g. in \cite{IA}, \cite{bressoud}, \cite{Szablowski2010(1)}.
\end{proof}

Two special cases $q\allowbreak=\allowbreak0,1$ are treated in the following Remark.

\begin{remark}
Assertions i)-iv) and ix) of the Lemma \ref{Hip} are true also for
$q\allowbreak=\allowbreak1.$ This follows elementary properties of Hermite
polynomials exposed e.g. in \cite{Andrews1999}. Further we have
\begin{align*}
\varphi\left(  x,t|1\right)   &  =\exp\left(  xt-t^{2}/2\right)
,\varphi\left(  x,t|0\right)  \allowbreak=\allowbreak\frac{1}{\left(
1-xt+t^{2}\right)  },\\
\tau\left(  x,t|y,\rho,1\right)  \allowbreak &  =\allowbreak\exp\left(
t\left(  x-\rho y\right)  -t^{2}(1-\rho^{2})/2\right)  ;\tau\left(
x,t|y,\rho,0\right)  \allowbreak=\allowbreak\frac{\left(  1-\rho yt+\rho
^{2}t^{2}\right)  }{\left(  1-xt+t^{2}\right)  }.
\end{align*}
For $q=1$ we have
\begin{equation}
f_{CN}\left(  x|y,\rho,1\right)  =\frac{1}{\sqrt{2\pi\left(  1-\rho
^{2}\right)  }}\exp\left(  -\frac{\left(  x-\rho y\right)  ^{2}}{2\left(
1-\rho^{2}\right)  }\right)  , \label{Norm2}%
\end{equation}
(i.e. Normal $N(\rho y,1-\rho^{2})$ distribution) while for $q\allowbreak
=\allowbreak0$ we have
\[
f_{CN}\left(  x|y,\rho,0\right)  \allowbreak=\allowbreak\frac{(1-\rho
^{2})\sqrt{4-x^{2}}}{2\pi((1-\rho^{2})^{2}-\rho(1+\rho^{2})xy+\rho^{2}%
(x^{2}+y^{2}))},
\]
$x,y\in\lbrack-2,2],$ $\left\vert \rho\right\vert <1,$ -the so called
Kesten--McKay distribution.
\end{remark}

As mentioned earlier we are going to consider continuous time generalization
of the process $\mathbf{X}$ considered in this section. Namely, more
precisely, we are going to prove the existence and present basic properties of
the process $\mathbf{Y\allowbreak=\allowbreak}\left\{  Y_{t}\right\}
_{t\in\mathbb{R}}$, satisfying the following condition:

\begin{condition}
[$\delta$]For every positive $\delta,$ random sequences
\begin{equation}
X_{n}^{\left(  \delta\right)  }=Y_{n\delta},~n\in\mathbb{Z},
\label{discreteversion}%
\end{equation}
are regular 1TSP.
\end{condition}

\section{Auxiliary properties of one dimensional random fields}

In this section we are going to make preparations for the proof of the
existence of process $\mathbf{Y}$. To do this we have to prove the following Lemma.

\begin{lemma}
\label{general}Let $\left\{  X_{i}\right\}  _{i\in\mathbb{Z}}$ 1TSP with
parameters $\rho$ and $q$. Let us fix $j\in\mathbb{N}$ and $m\in
\{0,\ldots,j-1\}$ and define $Z_{k}=X_{kj+m}$ for $k\in\mathbb{\mathbb{Z}}$.
Then $\left\{  Z_{k}\right\}  _{k\in\mathbb{Z}}$ is also 1TSP with parameters
$\rho^{j}$ and $q$.
\end{lemma}

To do this we need generalizations of properties (\ref{bryc1}) and
(\ref{bryc2}). They are given in Proposition presented below. On the other
hand this Proposition needs the following technical Lemma.

\begin{lemma}
\label{kwadraty}If $\mathbf{X}$ is a 1TSP then:
\begin{subequations}
\label{wzorki}%
\begin{align}
\mathbb{E}X_{n}^{4}  &  =\left(  2+q\right)  ,\label{4ty}\\
\mathbb{E}X_{n}^{2}X_{m}^{2}  &  =1+\rho^{2\left\vert n-m\right\vert }\left(
1+q\right)  ,\label{2_2}\\
\mathbb{E}X_{n}^{2}X_{n-j}X_{n+k}  &  =\rho^{j+k}\mathbb{E}X_{n}^{4}%
=\rho^{j+k}\left(  2+q\right)  . \label{2_1_1}%
\end{align}

\end{subequations}
\end{lemma}

\begin{proposition}
\label{u_bryc2}If $\mathbf{X}$ is a 1TSP then for $n,j,k\in\mathbb{\mathbb{N}%
}$

i) $\mathbb{E}\left(  X_{n}|\mathcal{F}_{\leq n-j,\geq n+k}\right)  $ is a
linear function of $X_{n-k}$ and $X_{n+j}$. More precisely we have
\begin{equation}
\mathbb{E}\left(  X_{n}|\mathcal{F}_{\leq n-k,\geq n+j}\right)  =\frac
{\rho^{j}\left(  1-\rho^{2k}\right)  }{1-\rho^{2\left(  j+k\right)  }}%
X_{n-j}+\frac{\rho^{k}\left(  1-\rho^{2j}\right)  }{1-\rho^{2\left(
j+k\right)  }}X_{n+k}. \label{liniowosc}%
\end{equation}

ii) $\mathbb{E}\left(  X_{n}^{2}|\mathcal{F}_{n-j\leq,\geq n+k}\right)  $ is a
linear function of $X_{n-j}^{2},$ $X_{n+k}^{2}$ and $X_{n-j}X_{n+k}$. In
particular%
\begin{equation}
\mathbb{E}\left(  X_{n}^{2}|\mathcal{F}_{\leq n-j,\geq n+k}\right)
=A_{jk}^{\left(  1\right)  }X_{n-j}^{2}+A_{jk}^{\left(  2\right)  }X_{n+k}%
^{2}+B_{jk}X_{n-j}X_{n+k}+C_{jk} \label{uog_bryc2}%
\end{equation}
where:
\begin{subequations}
\label{par_j}%
\begin{align}
A_{jk}^{\left(  1\right)  }  &  =\frac{\rho^{2j}\left(  1-\rho^{2k}\right)
\left(  1-q\rho^{2k}\right)  }{\left(  1-q\rho^{2\left(  j+k\right)  }\right)
\left(  1-\rho^{2\left(  j+k\right)  }\right)  },\label{_Aj1}\\
A_{jk}^{\left(  2\right)  }  &  =\frac{\rho^{2k}\left(  1-\rho^{2j}\right)
\left(  1-q\rho^{2j}\right)  }{\left(  1-q\rho^{2\left(  j+k\right)  }\right)
\left(  1-\rho^{2\left(  j+k\right)  }\right)  }\label{_Aj2}\\
B_{jk}  &  =\frac{\left(  q+1\right)  \rho^{\left(  j+k\right)  }\left(
1-\rho^{2j}\right)  \left(  1-\rho^{2k}\right)  }{\left(  1-q\rho^{2\left(
j+k\right)  }\right)  \left(  1-\rho^{2\left(  j+k\right)  }\right)
},\label{_Bj}\\
C_{jk}  &  =\frac{\left(  1-\rho^{2j}\right)  \left(  1-\rho^{2k}\right)
}{1-q\rho^{2\left(  j+k\right)  }}. \label{_Cj}%
\end{align}

\end{subequations}
\end{proposition}

Lengthy proofs of these facts as well as the proof of Lemma \ref{general} are
moved to section \ref{dowody}.

\section{$(\alpha,q)-$ Ornstein--Uhlenbeck and $q-$Wiener processes}

\subsection{Existence}

In this subsection we are going to prove the following Theorem

\begin{theorem}
\label{existence}The $L_{2}-$continuous\footnote{$L_{2}$ -continuous means
mean-square continuous} and stationary process $\mathbf{Y\allowbreak
=\allowbreak\allowbreak}\left\{  Y_{t}\right\}  _{t\in\mathbb{R}}$ that
satisfies for every $\delta>0$ Condition ($\delta$), exists. Moreover there
exist two numbers $q\in(-1,1]$ and $\alpha>0$ such that $\mathbf{Y}$ is Markov
with the marginals having density $f_{N}\left(  x|q\right)  $ and the
transition distribution having density $f_{CN}\left(  x|y,e^{-\alpha\left\vert
s-t\right\vert },q\right)  $ (i.e. $Y_{s}|Y_{t}=y\sim f_{CN}\left(
x|y,e^{-\alpha\left\vert s-t\right\vert },q\right)  $)
\end{theorem}

\begin{proof}
An easy but long proof is shifted to Section \ref{dowody}.
\end{proof}

In the sequel, when considering the continuous time generalizations of 1TSP we
will need the following generalization of 'non-singularity of covariance
matrix' assumption considered in the case of 1TSP: Let $\mathcal{X}=\left(
X_{t}\right)  _{t\in\mathbb{R}}$ be square integrable stochastic process and%
\begin{equation}
\forall n\in\mathbb{N};0\leq t_{1}<t_{2}<\ldots<t_{n}~\text{random variables
}X_{t_{1}},\ldots,X_{t_{n}}~\text{are linearly independent} \label{lin_ind}%
\end{equation}
which we will also refer to as linear independence assumption be satisfied by
$\mathcal{X}$.

\subsection{$(\alpha,q)-$OU processes}

Process $\mathbf{Y}$ with parameters $(\alpha,q)$ will be called continuous
time $(\alpha,q)-$ OU-process. (OU standing for Ornstein--Uhlenbeck). Let us
summarize what properties of $\mathbf{Y}$ can be deduced from the properties
of the discrete time regular 1TSP processes presented in \cite{bryc1},
\cite{bryc2}, \cite{matszab}, \cite{matszab2}, \cite{bms} and \cite{szab10}
(Corollary 6 p.13). Some of these properties also appeared in \cite{Bo} as a
by-product of considering some noncommutative model or in \cite{BryMaWe} in
quadratic harnesses context. Hence some of the properties of are not new but
they are scattered in the literature and we bring them together and collect in
groups devoted to particular features of analyzed processes.

As before let us define the following $\sigma-$fields defined by $\mathbf{Y}$.
$\mathcal{F}_{\leq s}:=\sigma\left(  X_{t}:t\leq s\right)  ,$ $\mathcal{F}%
_{\geq s}:=\sigma\left(  X_{t}:t\geq s\right)  ,$ and $\mathcal{F}_{\leq
s,\geq t}:=\sigma\left(  X_{\tau}:\tau\leq s\vee\tau\geq t\right)  $ for $s<t$.

Theorem below describes marginal and conditional distributions, presents
polynomials that are orthogonal with respect to these distributions as well as
gives conditional moments with respect to one-sided ($\mathcal{F}_{\leq s}$
and $\mathcal{F}_{\geq s})$ and two-sided ($\mathcal{F}_{\leq s,\geq t})$
$\sigma-$fields.

\begin{theorem}
\label{podstawowe}Let $\mathbf{Y}$ be a continuous time $(\alpha,q)-$
OU-process, $-1<q\leq1,$ $\alpha>0$. Then its state space is $S\left(
q\right)  $ and:

\begin{enumerate}
\item
\begin{align*}
\forall t  &  \in\mathbb{R}:Y_{t}\sim f_{N}\left(  x|q\right)  ,\\
\forall s  &  >t:Y_{s}|Y_{t}=y\sim f_{CN}\left(  x|y,e^{-\alpha\left(
s-t\right)  },q\right)  .
\end{align*}

\item $\mathbf{Y}$ is a stationary Markov process with $f_{CN}\left(
x|y,e^{-\alpha\left(  s-t\right)  },q\right)  $ as the density of its
transition probability.

\item $\mathbf{Y}$ is time symmetric. Moreover we have for any $n\in
\mathbb{N},s\in\mathbb{R}$ and $\delta,\gamma>0:$%
\begin{equation}
\mathbb{E}\left(  H_{n}\left(  Y_{s}\right)  |\mathcal{F}_{\leq s-\delta,\geq
s+\gamma}\right)  \allowbreak=\allowbreak\sum_{r=0}^{\left\lfloor
n/2\right\rfloor }\sum_{l=0}^{n-2r}A_{r,-\left\lfloor n/2\right\rfloor
+r+l}^{\left(  n\right)  }H_{l}\left(  Y_{s-\delta}|q\right)  H_{n-2r-l}%
\left(  Y_{s+\gamma}|q\right)  , \label{warmiedzy}%
\end{equation}
where $\left\lfloor \frac{n+2}{2}\right\rfloor \left\lfloor \frac{n+3}%
{2}\right\rfloor $ constants $A_{r,m}^{(n)};$. $r\allowbreak=\allowbreak
0,\ldots,\left\lfloor n/2\right\rfloor ,$ $m\allowbreak=\allowbreak
-\left\lfloor n/2\right\rfloor +r,\ldots,-\left\lfloor n/2\right\rfloor
+r\allowbreak+\allowbreak n-2r,$ depend only on $n,$ $q,e^{-\alpha\delta}$ and
$e^{-\alpha\gamma}.$\footnote{Since the paper was written and submitted the
exact form of coefficients $A_{r,s}^{\left(  n\right)  }$ can be derived from
the result presented in \cite{szab2011}, Thm. 2 following observation that the
conditional distribution of $X_{\sigma}$ given $X_{\sigma-\delta}$ ,
$X_{\sigma+\gamma}$ has Askey-Wilson density with specific complex
parameters.}

\item Families of polynomials $\left\{  H_{n}\left(  x|q\right)  \right\}
_{n\geq0}$ and $\left\{  P_{n}\left(  x|y,e^{-\alpha\left\vert s-t\right\vert
},q\right)  \right\}  _{n\geq0}$ given respectively by (\ref{Her}) and
(\ref{alsalamy}), are orthogonal polynomials of distributions defined by
(\ref{q_norm}) and (\ref{w_kowa}) respectively. That is in particular we have%
\begin{align}
\forall n  &  \geq1,t\in\mathbb{R}:\mathbb{E}H_{n}\left(  Y_{t}|q\right)
=0\label{_hermi}\\
\forall n  &  \geq1,s>t:\mathbb{E}\left(  P_{n}\left(  Y_{s}|Y_{t}%
,e^{-\alpha(s-t)},q\right)  |\mathcal{F}_{\leq t}\right)  =0\ \ a.s.~,
\label{_alsalam}%
\end{align}

\item
\begin{equation}
\forall n\geq1,s>t:\mathbb{E(}H_{n}\left(  Y_{s}|q\right)  |\mathcal{F}_{\leq
t})=e^{-n\alpha\left(  s-t\right)  }H_{n}\left(  Y_{t}|q\right)  ~~~a.s.~.
\label{herwar}%
\end{equation}

\item $\forall n\geq1:\left(  H_{n}\left(  Y_{t}|q\right)  \right)
_{t\in\mathbb{R}}$ is a stationary random process with covariance function
\[
K_{n}\left(  s,t\right)  =\left[  n\right]  _{q}!e^{-n\alpha\left\vert
s-t\right\vert }%
\]
and
\[
S_{n}\left(  \omega|\alpha\right)  =\frac{2n\alpha\left[  n\right]  _{q}%
!}{\omega^{2}+n^{2}\alpha^{2}},
\]
as its spectral density.
\end{enumerate}
\end{theorem}

\begin{proof}
[Remarks concerning the proof](1) and (2) are given in \cite{bryc1} and
\cite{bms}. (3) is given in \cite{szab10} (Corollary 5), (4) (\ref{_hermi}) is
given in \cite{IA} (but also in \cite{bryc1}) and (\ref{_alsalam}) is given in
\cite{bms}. (5) is given in \cite{bryc1} (6) Notice that from (5) it follows
that $\forall n,m\geq1,s,t\in\mathbb{R}\allowbreak:\allowbreak\mathbb{E}%
H_{m}\left(  Y_{t}|q\right)  H_{n}\left(  Y_{s}|q\right)  \allowbreak
=\allowbreak\delta_{\left\vert n-m\right\vert }e^{-n\alpha\left\vert
s-t\right\vert }\mathbb{E}H_{n}^{2}\left(  Y_{t}|q\right)  \allowbreak
=\allowbreak\delta_{\left\vert n-m\right\vert }e^{-n\alpha\left\vert
s-t\right\vert }\left[  n\right]  _{q}!,$ since $\mathbb{E}H_{n}^{2}\left(
Y_{t}|q\right)  \allowbreak=\allowbreak\left[  n\right]  _{q}!,$ by Lemma
\ref{Hip} (i). Further we use spectral decomposition theorem.
\end{proof}

\begin{remark}
As an immediate consequence of the assertion 1. of the above mentioned theorem
we see that if $\mathbf{Y}$ is a certain $\left(  \alpha,q\right)  -$OU
process then the process $\mathbf{Z}$ defined by re-scaling time in the
following way $Z_{s}\allowbreak=\allowbreak Y_{\alpha s},$ $s\in\mathbb{R}$ is
$\left(  1,q\right)  -$OU process.
\end{remark}

\begin{remark}
As it follows from Theorem \ref{podstawowe}, $4.$ the $\left(  \alpha
,q\right)  -$OU process is stationary and time homogenous. Thus its transition
operator is defined by
\[
P_{s,t}(f)\left(  y\right)  =\int_{S\left(  q\right)  }f\left(  x\right)
f_{CN}\left(  x|y,e^{-\alpha(t-s)},q\right)  dx
\]
for a function $f$ defined below and $t>s$ depends in fact on $t-s=\tau.$ Let
us define
\[
P^{\tau}\left(  .\right)  \allowbreak=\allowbreak P_{t,t+\tau}\left(
.\right)  .
\]
Operators $\left\{  P^{\tau}\right\}  _{\tau>0}$ form a semigroup of operators
as it follows from Lemma \ref{Hip}, iv.
\end{remark}

As conclusions of assertions of Theorem \ref{podstawowe} we have the following
Theorem that contains their implications to the Markovian properties of the
process $\mathbf{Y}$\textbf{.}

Before we formulate appropriate theorem we need to introduce some additional notation.

Let $B(q)\allowbreak=\allowbreak L_{2}\left(  S\left(  q\right)
,\mathcal{B}\left(  S\left(  q\right)  \right)  ,P_{N}\left(  q\right)
\right)  ,$ where $\mathcal{B}\left(  S\left(  q\right)  \right)  $ denotes
$\sigma-$field of Borel subsets of $S\left(  q\right)  $ and $P_{N}\left(
q\right)  $ denotes measure with density $f_{N}\left(  .|q\right)  $.
Obviously we have
\[
B\left(  q\right)  \allowbreak=\allowbreak\left\{  f:S\left(  q\right)
\longrightarrow S\left(  q\right)  :f\left(  x\right)  \allowbreak
=\allowbreak\sum_{j=0}^{\infty}\frac{b_{j}}{\sqrt{\left[  i\right]  _{q}!}%
}H_{j}\left(  x|q\right)  ;\sum_{j\geq0}\left\vert b_{j}\right\vert
^{2}<\infty\right\}  .
\]

Let us further denote:
\[
B^{0}\left(  q\right)  \allowbreak=\allowbreak\{f\in B\left(  q\right)
\allowbreak:\allowbreak f\left(  x\right)  \allowbreak=\allowbreak\sum
_{j=0}^{\infty}\frac{b_{j}}{\sqrt{\left[  j\right]  _{q}!}}H_{j}\left(
x|q\right)  ;\sum_{j\geq j}j^{2}\left\vert b_{j}\right\vert ^{2}<\infty\}.
\]
We have $B^{0}\left(  q\right)  \subset B\left(  q\right)  $ .

Let $A$ denote infinitesimal operator of the $\left(  \alpha,q\right)  -$OU process.

\begin{theorem}
\label{ou-_wlasn}Let $\mathbf{Y}$ be $\left(  \alpha,q\right)  -$OU process
with $\left\vert q\right\vert <1$. Then:

i) $\operatorname*{var}\left(  Y_{t}\right)  \allowbreak=\allowbreak1,$ and
for $s\geq0,$ $\operatorname*{var}\left(  Y_{t+s}|Y_{t}\right)  \allowbreak
=\allowbreak1-e^{-2\alpha s}$. Hence in particular trajectories of $\left(
\alpha,q\right)  -$OU process are c\`{a}dl\`{a}g functions with values in
$S\left(  q\right)  $,

ii) transition operator $P^{\tau}$ for $\tau>0$ is defined by the following
relationship
\begin{equation}
B\left(  q\right)  \ni f\left(  x\right)  \allowbreak=\allowbreak\sum_{j\geq
0}\frac{b_{j}}{\sqrt{\left[  j\right]  _{q}!}}H_{j}\left(  x|q\right)
\allowbreak->\allowbreak P^{\tau}(f)\left(  x\right)  \allowbreak
=\allowbreak\sum_{j\geq0}e^{-\alpha\tau j}\frac{b_{j}}{\sqrt{\left[  j\right]
_{q}!}}H_{j}\left(  x|q\right)  , \label{formula}%
\end{equation}

iii) the family of transitional probabilities $\left\{  P^{\tau}\right\}
_{\tau>0}$ is Feller -continuous. In particular process $\mathbf{Y}$ has
strong Markov property,

iv) the family $\left\{  P^{\tau}\right\}  _{\tau>0}$ is right continuous
consequently process $\mathbf{Y}$ is a Feller process. Moreover its
infinitesimal operator $A$ exists and is defined on a subset $B^{0}$ by the
following formula:
\[
B^{0}\left(  q\right)  \allowbreak\ni\allowbreak f\left(  x\right)
\allowbreak=\allowbreak\sum_{j\geq0}\frac{b_{j}}{\sqrt{\left[  j\right]
_{q}!}}H_{j}\left(  x|q\right)  \allowbreak->\allowbreak A\left(  f\right)
\left(  x\right)  \allowbreak=\allowbreak-\alpha\sum_{j\geq1}\frac{jb_{j}%
}{\sqrt{\left[  j\right]  _{q}!}}H_{j}\left(  x|q\right)  \allowbreak
\in\allowbreak B\left(  q\right)  .
\]

\end{theorem}

\begin{proof}
The proof is shifted to section \ref{dowody}.
\end{proof}

\begin{remark}
Let us remark that recently Anshelevich in \cite{Ansh2011} (Lemma 20)
expressed the infinitesimal operator in an integral form. To be precise he
expressed in this form the infinitesimal operator of the so called $q-$Wiener
($q-$Brownian motion as he calls it) process to be considered in the next
subsection and related to $\left(  \alpha,q\right)  -$OU by the continuous
transformation (\ref{def_q-W}).
\end{remark}

The detailed forms of constants $A_{r,m}^{\left(  n\right)  }$ defined by
(\ref{herwar}) are given below following Corollary 6 of \cite{szab10}.

\begin{corollary}
\label{warposredni}$A_{0,-\left\lfloor n/2\right\rfloor +l}^{(n)}%
\allowbreak=\allowbreak%
\genfrac{[}{]}{0pt}{}{n}{l}%
_{q}\frac{e^{-(n-l)\alpha\delta}\left(  e^{-2\alpha\gamma}\right)
_{n-l}e^{-l\alpha\gamma}\left(  e^{-2\alpha\delta}\right)  _{l}}{\left(
e^{-2\alpha\left(  \delta+\gamma\right)  }\right)  _{n}},$ $l\allowbreak
=\allowbreak0,\ldots,n,$ $n\allowbreak=\allowbreak1\ldots,4$. If
$n\allowbreak\leq3$ then $A_{1,-\left\lfloor n/2\right\rfloor +l}^{\left(
n\right)  }\allowbreak=\allowbreak-\left[  n-1\right]  _{q}e^{-\alpha\left(
\delta+\gamma\right)  }A_{0,-\left\lfloor n/2\right\rfloor +l}^{\left(
n\right)  },$ $l\allowbreak=\allowbreak1,\ldots,n-1,$ If $n\allowbreak
=\allowbreak4$ then $A_{1,j}^{\left(  4\right)  }\allowbreak=\allowbreak
-\left[  3\right]  _{q}e^{-\alpha\left(  \delta+\gamma\right)  }A_{0,j}%
^{(4)},$ $j\allowbreak=\allowbreak-1,1$ and $A_{1,0}^{(4)}\allowbreak
=\allowbreak-\left[  2\right]  _{qq}^{2}e^{-\alpha\left(  \delta
+\gamma\right)  }A_{0,0}^{(4)},$ $A_{2,0}^{(4)}\allowbreak=\allowbreak
q(1+q)_{q}e^{-2\alpha\left(  \delta+\gamma\right)  }A_{0,0}^{(4)}.$
\end{corollary}

In particular we obtain known (see e.g. \cite{bryc05}) formula
\begin{gather}
\operatorname{var}\left(  Y_{s}|\mathcal{F}_{\leq s-\delta,\geq s+\gamma
}\right)  \allowbreak=\allowbreak\frac{(1-e^{-2\alpha\delta})\left(
1-e^{-2\alpha\gamma}\right)  }{\left(  1-qe^{-2\alpha(\delta+\gamma)}\right)
}\allowbreak\label{warvar1}\\
\times\allowbreak(1-\frac{\left(  1-q\right)  (X_{i-1}-e^{-\alpha
(\delta+\gamma)}X_{i+1})\left(  X_{i}-e^{-\alpha(\delta+\gamma)}%
X_{i-1}\right)  }{\left(  1-e^{-2\alpha(\delta+\gamma)}\right)  ^{2}}).
\label{warvar2}%
\end{gather}

Following assertion vii of Lemma \ref{Hip} and we can express assertion
Theorem \ref{podstawowe}, (5) in the following martingale-like form.

\begin{theorem}
\label{martyngaly}For $|q|<1$ we have%
\[
\forall\gamma^{2}\left(  1-q\right)  <1,s>t:\mathbb{E}\left(  \varphi\left(
Y_{s},\gamma|q\right)  |\mathcal{F}_{\leq t}\right)  \allowbreak
=\allowbreak\varphi\left(  Y_{t},\gamma e^{-\alpha(s-t)}|q\right)
~\text{a.s.}%
\]
where positive function $\varphi$ is defined in Lemma \ref{Hip}, vii).

In particular for $q\allowbreak=\allowbreak0$ (so called free $\alpha
-$OU-process) we have
\[
\forall\left\vert \gamma\right\vert <1,s>t:\mathbb{E}\left(  \frac{1}{1-\gamma
Y_{s}+\gamma^{2}}|\mathcal{F}_{\leq t}\right)  \allowbreak=\allowbreak\frac
{1}{1-\gamma e^{-\alpha(t-s)}Y_{t}+\gamma^{2}e^{-2\alpha\left(  t-s\right)  }%
}~~a.s.,
\]
while for $q\allowbreak=\allowbreak1$ we get well known formula:%
\begin{align*}
\forall\gamma &  \in\mathbb{R},s>t:\mathbb{E}\left(  \exp\left(  \gamma
Y_{s}-\gamma^{2}/2\right)  |\mathcal{F}_{\leq t}\right) \\
&  =\exp\left(  \gamma e^{-\alpha\left(  s-t\right)  }Y_{t}-\left(  \gamma
e^{-\alpha\left(  s-t\right)  }\right)  ^{2}/2\right)  ~a.s.
\end{align*}

\end{theorem}

Strict proof is very much alike the proof of Corollary \ref{rozk_q_wienera}
below, so we will not present it here.

\begin{remark}
\label{zalozenia_q-wienera}Notice that following considerations of the section
\ref{wlasnosci}, concerning existence of $(\alpha,q)-$OU process one needed
only the following \emph{two,} symmetric in time, conditions apart from linear independence:

\begin{enumerate}
\item $\forall n\in\mathbb{Z},d>0:$%
\[
\mathbb{E}\left(  Y_{nd}|\mathcal{F}_{\leq\left(  n-1\right)  d,\geq\left(
n+1\right)  d}\right)  =\frac{e^{-\alpha d}}{1+e^{-2\alpha d}}\left(
Y_{\left(  n-1\right)  d}+Y_{\left(  n+1\right)  d}\right)  ,
\]

\item $\forall n\in\mathbb{Z},d>0:$%
\[
\mathbb{E}\left(  Y_{nd}^{2}|\mathcal{F}_{\leq\left(  n-1\right)
d,\geq\left(  n+1\right)  d}\right)  =\hat{A}\left(  Y_{\left(  n-1\right)
d}^{2}+Y_{\left(  n+1\right)  d}^{2}\right)  +\hat{B}Y_{\left(  n-1\right)
d}Y_{\left(  n+1\right)  d}+\hat{C},
\]
where $\hat{A}=A_{1}\left(  d,d\right)  =A_{2}\left(  d,d\right)  ,$ $\hat
{B}=B\left(  d,d\right)  ,$ $\hat{C}=C\left(  d,d\right)  $. That is we need
only symmetric (and discrete for all increments $d>0)$ versions of condition
defining $\mathbb{E}\left(  Y_{s}^{2}|\mathcal{F}_{\leq s-\delta,\geq
s+\gamma}\right)  $.
\end{enumerate}
\end{remark}

\begin{remark}
\label{no-existance}Notice that continuous time $\left(  \alpha,q\right)  -$OU
process does not exist for $q>1$. It is so because for discrete time 1TSP with
parameters $(q,\rho)$ the following relationship between parameters $q$ and
$\rho$ must be satisfied for some integer $n$ : $\rho^{2}q^{n}=1$ or
equivalently ratio $\frac{\log q}{\log\rho^{2}}$ must be equal to some
integer. However if $\left(  \alpha,q\right)  -$OU process existed parameter
$\rho$ would depend on time parameter $t$ in the following way $\rho^{2}%
=\exp\left(  -2\alpha t\right)  $ for some fixed positive $\alpha$ and
consequently $\frac{\log q}{\alpha t}$ would have to be integer for all real
$t$ which is impossible.
\end{remark}

\subsection{$q-$Wiener process}

Let $\mathbf{Y}$ be given $(\alpha,q)-$OU process. Let us define:%
\begin{equation}
X_{0}=0;~\forall\tau>0:X_{\tau}=\sqrt{\tau}Y_{\log\tau/2\alpha}.
\label{def_q-W}%
\end{equation}
Process $\mathbf{X=}\left(  X_{\tau}\right)  _{\tau\geq0}$ will be called
$q-$Wiener process. Let us also introduce the following filtration:
\[
\mathcal{F}_{\leq\theta}^{X}=\sigma\left(  X_{\tau}:\tau\leq\theta\right)
\allowbreak=\allowbreak\sigma\left(  Y_{t}:t\leq\log\theta/2\alpha\right)
\allowbreak\left(  =\mathcal{F}_{\leq\log\theta/2\alpha}\right)  .
\]

\begin{remark}
From (\ref{def_q-W}) and from the properties of $(\alpha,q)-$OU process it
follows that $\mathbf{X}$ is a self-similar process since for $c>0$ we see
that $\left\{  Z_{\tau}\right\}  _{\tau\geq0}$ where $Z_{\tau}\allowbreak
=\allowbreak c^{-1}X_{c^{2}\tau};$ $\tau\geq0,$ is also a $q-$Wiener process.
\end{remark}

Following definition given by (\ref{def_q-W}) and Theorem \ref{podstawowe} of
the previous section we have the following Theorem whose detailed proof is
simple but lengthy and thus is shifted to section \ref{dowody}:

\begin{theorem}
\label{wl_q-Wienera}Let $\mathbf{X}$ be a $q-$Wiener process, then it has the
following properties:

\begin{enumerate}
\item $\forall\tau,\sigma\geq0:\operatorname{cov}\left(  X_{\tau},X_{\sigma
}\right)  =\min\left(  \tau,\sigma\right)  $.

\item $\forall\tau>0:X_{\tau}\sim\frac{1}{\sqrt{\tau}}f_{N}\left(  \frac
{x}{\sqrt{\tau}}|q\right)  $.

\item For $\tau>\sigma:X_{\tau}-X_{\sigma}|X_{\sigma}=y\sim\frac{1}{\sqrt
{\tau}}f_{CN}\left(  \frac{x+y}{\sqrt{\tau}}|\frac{y}{\sqrt{\sigma}}%
,\sqrt{\frac{\sigma}{\tau}},q\right)  $.

\item For all $n\geq1$ and $0<\sigma\leq\tau$ we have
\begin{align}
\mathbb{E}\left(  \tau^{n/2}H_{n}\left(  \frac{X_{\tau}}{\sqrt{\tau}%
}|q\right)  |\mathcal{F}_{\leq\sigma}^{X}\right)   &  =\sigma^{n/2}%
H_{n}\left(  \frac{X_{\sigma}}{\sqrt{\sigma}}|q\right)  ,a.s.~,\label{prop1}\\
\mathbb{E}\left(  \sigma^{-n/2}H_{n}\left(  \frac{X_{\sigma}}{\sqrt{\sigma}%
}|q\right)  \mathcal{F}_{\geq\tau}^{X}\right)   &  =\tau^{-n/2}H_{n}\left(
\frac{X_{\tau}}{\sqrt{\tau}}|q\right)  ~a.s.~. \label{prop2}%
\end{align}
Hence $\forall n\geq1$ the pair $\left(  Z_{\tau}^{\left(  n\right)
},\mathcal{F}_{\leq\tau}^{X}\right)  _{\tau\geq0},$ where $Z_{\tau}^{\left(
n\right)  }\allowbreak=\allowbreak\tau^{n/2}H_{n}\left(  X_{\tau}/\sqrt{\tau
}|q\right)  ,$ $\tau\geq0$ is a martingale and the pair $\left(  V_{\tau
}^{\left(  n\right)  },\mathcal{F}_{\geq\tau}^{X}\right)  _{\tau\geq0},$ where
$V_{\tau}^{\left(  n\right)  }\allowbreak=$ $\allowbreak\tau^{-n/2}%
H_{n}\left(  \frac{X_{\tau}}{\sqrt{\tau}}|q\right)  $ is the reverse
martingale. In particular $\mathbf{X}$ is a martingale and $\left(
\frac{X_{\tau}}{\tau},\mathcal{F}_{\geq\tau}^{X}\right)  _{\tau>0}$ is the
reversed martingale.

\item We have also $\forall n\geq1,$ $\delta,\gamma\geq0,$ $\sigma\geq\delta$%
\[
\mathbb{E}\left(  H_{n}\left(  \frac{X_{\sigma}}{\sqrt{\sigma}}|q\right)
|\mathcal{F}_{\leq\sigma-\delta,\geq\sigma+\gamma}\right)  \allowbreak
=\allowbreak\sum_{r=0}^{\left\lfloor n/2\right\rfloor }\sum_{l=0}%
^{n-2r}A_{r,-\left\lfloor n/2\right\rfloor +r+l}^{\left(  n\right)  }%
H_{l}\left(  \frac{X_{\sigma-\delta}}{\sqrt{\sigma-\delta}}|q\right)
H_{n-2r-l}\left(  \frac{X_{\sigma+\gamma}}{\sqrt{\sigma+\gamma}}|q\right)  ,
\]
where $\left\lfloor \frac{n+2}{2}\right\rfloor \left\lfloor \frac{n+3}%
{2}\right\rfloor $ constants $A_{r,s}^{(n)};$. $r\allowbreak=\allowbreak
0,\ldots,\left\lfloor n/2\right\rfloor ,$ $s\allowbreak=\allowbreak
-\left\lfloor n/2\right\rfloor +r,\ldots,-\left\lfloor n/2\right\rfloor
+r\allowbreak+\allowbreak n-2r$ depend only on $n,$ $q,$and numbers
$\sigma,\delta$ and $\gamma.$\footnote{See footnote to assertion 3 of Theorem
\ref{podstawowe}.}
\end{enumerate}
\end{theorem}

We have the following immediate, easy observation.

\begin{remark}
From assertion $2$ of the above mentioned Theorem we see that for every $t>0,$
$tX_{1/t}$ and $X_{t}$ have the same distribution.
\end{remark}

\begin{corollary}
\label{har4}Let $\mathbf{X}$ be a $q-$Wiener process . We have:

i)
\begin{align}
\mathbb{\forall\tau}  &  >\sigma>0:\mathbb{E}\left(  \left(  X_{\tau
}-X_{\sigma}\right)  ^{2}|\mathcal{F}_{\leq\sigma}^{X}\right)  =\tau
-\sigma~a.s.~,\label{2gi_wkowy}\\
\mathbb{\forall\tau}  &  >\sigma>0:\mathbb{E}\left(  \left(  X_{\tau
}-X_{\sigma}\right)  ^{3}|\mathcal{F}_{\leq\sigma}^{X}\right)  =-(1-q)\left(
\tau-\sigma\right)  X_{\sigma}~a.s.~,\label{3ci_wkowy}\\
\mathbb{\forall\tau}  &  >\sigma>0:\mathbb{E}\left(  \left(  X_{\tau
}-X_{\sigma}\right)  ^{4}|\mathcal{F}_{\leq\sigma}^{X}\right)
=\label{4ty_wkowy}\\
&  \left(  \tau-\sigma\right)  \left(  X_{\sigma}^{2}\left(  1-q\right)
^{2}+\left(  2+q\right)  \left(  \tau-\sigma\right)  +\sigma\left(
1-q^{2}\right)  \right)  ~a.s.~.\nonumber
\end{align}
Hence $<X_{\tau}>\allowbreak=\allowbreak\tau$ moreover $\mathbf{X}$ does not
have independent increments.

ii) Almost every path of process $\mathbf{X}$ has at any point left and right
hand side limits, thus can be modified to have c\`{a}dl\`{a}g trajectories.

iii) Process $\mathbf{X}$ is Feller -continuous process and has strong Markov property.
\end{corollary}

The following corollary gives more detailed consequences of \cite{szab10}
(Corollary 6 p.13). These conditional moments are known of course (see e.g.
\cite{BryMaWe} for $n\allowbreak=\allowbreak1,2$ since $q-$Wiener process, as
pointed out in the introduction can be obtained as a particular quadratic harness.

\begin{corollary}
\label{warobustr}We have $A_{0,-\left\lfloor n/2\right\rfloor +l}%
^{(n)}\allowbreak=\allowbreak%
\genfrac{[}{]}{0pt}{}{n}{l}%
_{q}\frac{\rho_{1}^{n-l}\left(  \rho_{2}^{2}\right)  _{n-l}\rho_{2}^{l}\left(
\rho_{1}^{2}\right)  _{l}}{\left(  \rho_{1}^{2}\rho_{2}^{2}\right)  _{n}},$
$l\allowbreak=\allowbreak0,\ldots,n,$ $n\allowbreak=\allowbreak1\ldots,4$. If
$n\allowbreak\leq3$ then $A_{1,-\left\lfloor n/2\right\rfloor +l}^{\left(
n\right)  }\allowbreak=\allowbreak-\left[  n-1\right]  _{q}\rho_{1}\rho
_{2}A_{0,-\left\lfloor n/2\right\rfloor +l}^{\left(  n\right)  },$
$l\allowbreak=\allowbreak1,\ldots,n-1,$ If $n\allowbreak=\allowbreak4$ then
$A_{1,j}^{\left(  4\right)  }\allowbreak=\allowbreak-\left[  3\right]
_{q}\rho_{1}\rho_{2}A_{0,j}^{(4)},$ $j\allowbreak=\allowbreak-1,1$ and
$A_{1,0}^{(4)}\allowbreak=\allowbreak-\left[  2\right]  _{q}^{2}\rho_{1}%
\rho_{2}A_{0,0}^{(4)},$ $A_{2,0}^{(4)}\allowbreak=\allowbreak q(1+q)\rho
_{1}^{2}\rho_{2}^{2}A,$ where $\rho_{1}\allowbreak=\allowbreak\sqrt
{(\sigma-\delta)/\sigma}$ , $\rho_{2}\allowbreak=\allowbreak\sqrt
{\sigma/\left(  \sigma+\gamma\right)  }.$ In particular we can deduce the
following known (from say \cite{Bo}, \cite{BryWe} or \cite{bryc05}) formulae:
\begin{subequations}
\label{harnesy}%
\begin{gather}
\mathbb{E}\left(  X_{\sigma}|\mathcal{F}_{\leq\sigma-\delta,\geq\sigma+\gamma
}\right)  =\frac{\gamma}{\delta+\gamma}X_{\sigma-\delta}+\frac{\delta}%
{\delta+\gamma}X_{\sigma+\gamma}\ \ a.s.~,\label{linear}\\
\mathbb{E}\left(  X_{\sigma}^{2}|\mathcal{F}_{\leq\sigma-\delta,\geq
\sigma+\gamma}\right)  =\frac{\delta\gamma}{\left(  \delta+\gamma\right)
\left(  \sigma\left(  1-q\right)  +\gamma+q\delta\right)  }\times
\label{dual1}\\
\left(  \left(  \left(  1-q\right)  \sigma+\gamma\right)  \frac{X_{\sigma
-\delta}^{2}}{\delta}+\left(  \left(  1-q\right)  \sigma+q\delta\right)
\frac{X_{\sigma+\gamma}^{2}}{\gamma}+\left(  q+1\right)  X_{\sigma-\delta
}X_{\sigma+\gamma}+\left(  \delta+\gamma\right)  \right)  . \label{dual2}%
\end{gather}
\end{subequations}
\begin{align}
\operatorname{var}\left(  X_{\sigma}|\mathcal{F}_{\leq\sigma-\delta,\geq
\sigma+\gamma}\right)   &  =\frac{\delta\gamma}{\left(  \sigma\left(
1-q\right)  +\gamma+q\delta\right)  }\label{warboth}\\
&  \times\left(  1-\left(  1-q\right)  \frac{\left(  X_{\sigma+\gamma
}-X_{\sigma-\delta}\right)  }{\delta+\gamma}\frac{\left(  \left(
\sigma+\gamma\right)  X_{\sigma-\delta}-\left(  \sigma-\delta\right)
X_{\sigma+\gamma}\right)  }{\delta+\gamma}\right)  .\nonumber
\end{align}

\end{corollary}

Following Theorem \ref{wl_q-Wienera} we get:

\begin{corollary}
\label{rozk_q_wienera} Let $\mathbf{X=}\left(  X_{\tau}\right)  _{\tau\geq0}$
be a $q-$Wiener process. $\forall$ $s\in\mathbb{R},$ $0<\sigma^{-1},\tau
<\frac{1}{s^{2}(1-q)},$ the following pairs:%
\begin{align}
&  \left(  \varphi\left(  \frac{X_{\tau}}{\sqrt{\tau}},s\sqrt{\tau}|q\right)
,\mathcal{F}_{\leq\tau}^{X}\right)  _{1/((1-q)s^{2})>\tau\geq0}%
,\label{general_martyngal}\\
&  \left(  \varphi\left(  \frac{X_{\sigma}}{\sqrt{\sigma}},\frac{s}%
{\sqrt{\sigma}}|q\right)  ,\mathcal{F}_{\geq\sigma}^{X}\right)  _{\sigma
>(1-q)s^{2}\geq0}. \label{inverse_martyngal}%
\end{align}
are positive respectively martingale and reversed martingale, where function
$\varphi(x,t|q)$ is a characteristic function of $q-$Hermite polynomials and
is defined in Lemma \ref{Hip}, vii).

In particular we get:
\begin{align*}
\forall s  &  \in\mathbb{R},\frac{1}{(1-q)s^{2}}>\tau>0:\mathbb{E}\left(
\varphi\left(  \frac{X_{\tau}}{\sqrt{\tau}},s\sqrt{\tau}|q\right)  \right)
=1,\\
\forall s  &  \in\mathbb{R},s^{2}(1-q)<\tau:\mathbb{E}\left(  \varphi\left(
\frac{X_{\tau}}{\sqrt{\tau}},\frac{s}{\sqrt{\tau}}|q\right)  \right)  =1.
\end{align*}
$0-$Wiener process (sometimes called free Wiener process) satisfies:%
\begin{align*}
\forall s  &  \in\mathbb{R}\mathbf{,}0<\tau<\sigma<\frac{1}{s^{2}}%
:\mathbb{E}\left(  \frac{1}{1-sX_{\tau}+\tau s^{2}}|\mathcal{F}_{\leq\sigma
}^{X}\right)  =\frac{1}{1-sX_{\sigma}+\sigma s^{2}};~a.s.\\
\forall s  &  \in\mathbb{R}\mathbf{,}0<s^{2}<\sigma<\tau:\mathbb{E}\left(
\frac{1}{1-sX_{\sigma}/\sigma+s^{2}/\sigma}|\mathcal{F}_{\geq\tau}^{X}\right)
=\frac{1}{1-sX_{\tau}/\tau+s^{2}/\tau};~a.s.
\end{align*}

$1-$Wiener process satisfies of course%
\begin{align*}
\forall s  &  \in\mathbb{R}:\mathbb{E}\left(  \exp\left(  sX_{\tau}-\tau
s^{2}/2\right)  |\mathcal{F}_{\leq\sigma}^{X}\right)  =\exp\left(  sX_{\sigma
}-\sigma s^{2}/2\right)  ;~a.s.\\
\forall s  &  \in\mathbb{R}:\mathbb{E}\left(  \exp\left(  sX_{\sigma}%
/\sigma-s^{2}/\left(  2\sigma\right)  \right)  |\mathcal{F}_{\geq\tau}%
^{X}\right)  =\exp\left(  sX_{\tau}/\tau-s^{2}/\left(  2\tau\right)  \right)
;~a.s.
\end{align*}

\end{corollary}

\begin{remark}
Let us recall following \cite{BryMaWe} that quadratic harnesses are roughly
speaking square integrable processes $\left\{  Z_{t}\right\}  _{t\in
\mathbb{R}^{+}}$ that satisfy 5 conditions imposed on their covariance
functions and on the first two conditional moments given the past and given
the past and the future. These conditions require that $\mathbb{E}\left(
Z_{t}|\mathcal{F}_{\leq t-\delta}\right)  $ and $\mathbb{E(}Z_{t}%
|\mathcal{F}_{\leq t-\delta,\geq t+\gamma})$ be linear functions while
$\mathbb{E}\left(  Z_{t}^{2}|\mathcal{F}_{\leq t-\delta}\right)  $ and
$\mathbb{E(}Z_{t}^{2}|\mathcal{F}_{\leq t-\delta,\geq t+\gamma})$ be a
quadratic functions of respectively $Z_{t-\delta}$ and $Z_{t-\delta}$ and
$Z_{t+\gamma}.$ Examining assertion of Corollary \ref{warobustr} we see that
$q-$Wiener process has this property. Thus we deduce that $\mathbf{X}$ is a
quadratic harness with parameters (introduced in \cite{BryMaWe})
$\theta\allowbreak=\allowbreak\eta\allowbreak=\allowbreak\tau\allowbreak
=\allowbreak\sigma=0$.
\end{remark}

\begin{remark}
From Theorem \ref{wl_q-Wienera}, 4. it follows that $\left(  X_{\tau
},\mathcal{F}_{\leq\tau}^{X}\right)  _{\tau\geq0},$ $\left(  X_{\tau}^{2}%
-\tau,\mathcal{F}_{\leq\tau}^{X}\right)  _{\geq0}$ are martingales and
$\left(  X_{\tau}/\tau,\mathcal{F}_{\geq\tau}^{X}\right)  _{\tau>0}$ $\left(
X_{\tau}^{2}/\tau^{2}-1/\tau,\mathcal{F}_{\geq\tau}^{X}\right)  _{\tau>0}$ are
reversed martingales. Thus if main result of Weso\l owski's paper \cite{Wes}
(stating that a process satisfying these conditions should be a Wiener
process) was true we would deduce that $\mathbf{X}$ is the Wiener process. But
it is not at least for $|q|<1.$ Let us note that Weso\l owski is aware of this
since in \cite{bryc05} he (together with Bryc) gives examples that contradict
the result of \cite{Wes}. He however did not point out were was the mistake.
Recall that Weso\l owski considered processes that have the property that they
and their squares were both martingales and reversed martingales.
Weso\l owski's argument was based on the value $\mathbb{E}\left(  X_{\tau
}-X_{\sigma}\right)  ^{4}$ that had to be calculated for the considered
process. Following formulae presented above we deduce that $\mathbb{E}\left(
X_{\tau}-X_{\sigma}\right)  ^{4}\allowbreak=\allowbreak\left(  2+q\right)
\left(  \tau-\sigma\right)  ^{2}+2(1-q)\sigma\left(  \tau-\sigma\right)  $ for
a $q-$Wiener process, while Weso\l owski stated that $\mathbb{E}\left(
X_{\tau}-X_{\sigma}\right)  ^{4}\allowbreak=\allowbreak c\left(  \tau
-\sigma\right)  ^{2}+2\left(  c-3\right)  \sigma\left(  \sigma+\tau\right)  ,$
where $c\allowbreak=\allowbreak\mathbb{E}X_{\tau}^{4}/\tau^{2}.$ In the case
of $q-$Wiener process one can calculate that $c\allowbreak=\allowbreak(2+q)$
thus Weso\l owski made mistake in his calculations.
\end{remark}

\begin{remark}
\label{kowariancje}Following (\ref{2gi_wkowy}) and (\ref{3ci_wkowy}) we have
\begin{gather*}
\forall0<\sigma<\tau\leq\upsilon<\omega:\operatorname*{cov}\left(  X_{\tau
}-X_{\sigma},X_{\omega}-X_{\upsilon}\right)  =0,\\
\operatorname*{cov}\left(  \left(  X_{\tau}-X_{\sigma}\right)  ^{2},\left(
X_{\omega}-X_{\upsilon}\right)  ^{2}\right)  =0,\\
\operatorname{cov}\left(  \left(  X_{\tau}-X_{\sigma}\right)  ^{3},\left(
X_{\omega}-X_{\upsilon}\right)  ^{3}\right)  =-\left(  1-q\right)  \left(
\tau-\sigma\right)  \left(  \omega-\upsilon\right)  \left(  \tau\left(
2+q\right)  -\sigma\left(  1+2q\right)  \right)  .
\end{gather*}

\end{remark}

\section{Open problems}

As it follows from the description of $q-$Wiener and $(\alpha,q)-$ OU
processes they do not allow continuous paths modifications. Their paths have
jumps. Besides both left- and right-hand side limits exist at any jumping
point. Consequently the paths of these processes do not have discontinuities
of the second kind. Thus there are immediate several questions:

\begin{enumerate}
\item In general on every finite interval there can be infinitely many jumps.
Is it true? Or can one prove some additional properties of these processes
that would eliminate this case? Certainly such properties do not exist for all
$\left\vert q\right\vert <1.$ The case $q\allowbreak=\allowbreak0$ leads to
the Cauchy process that has infinitely many jumps on every finite interval.
But may be one can find $q_{0}$ such that for $q_{0}<q<1$ the $q-$Wiener
process has only finite number of jumps on every finite interval?

\item What is the distribution of the size of jumps of $(\alpha,q)-$ OU
process. It is stationary. But is it of continuous, discrete or singular type.
Or may be it is a mixture?

\item On the other hand for every $\varepsilon>0$ there are only finite number
of jumps of the size not less than $\varepsilon.$ It follows from the symmetry
$(\alpha,q)-$ OU processes with respect to time argument that inter jumps
intervals between such jumps have the same distributions. What is the
\emph{distribution of the length} of those intervals. Strong Markov property
would suggest exponential distribution. Is it true? Do they form a
\emph{renewal process} i.e. are those intervals independent? Probably not, but
it needs justification.

\item What are the properties of quadratic variations of martingales
associated with $q$-Wiener processes.

\item Recently Anshelevich et al. (see \cite{ANBo}) have proved so called free
infinite divisibility of $q-$Normal distribution. Is transitional distribution
with density $f_{CN}$ also free infinitely divisible?
\end{enumerate}

\section{Proofs of the results\label{dowody}}

\begin{proof}
[Proof of Lemma \ref{kwadraty}]1. Remembering that $H_{4}\left(  x|q\right)
\allowbreak=\allowbreak x^{4}\allowbreak-\allowbreak(3+2q+q^{2})x^{2}%
\allowbreak+\allowbreak(1+q+q^{2})$ we have
\[
x^{4}=H_{4}\left(  x|q\right)  +(3+2q+q^{2})H_{2}\left(  x|q\right)  +2+q.
\]
Hence $\mathbb{E}X_{n}^{4}=2+q$ since $\mathbb{E}H_{4}\left(  X_{m}\right)
\allowbreak=\allowbreak\mathbb{E}H_{2}\left(  X_{m}\right)  \allowbreak
=\allowbreak0$. Now using the fact that $\mathbb{E}\left(  H_{2}\left(
X_{n}\right)  |\mathcal{F}_{\leq m}\right)  =\rho^{2(n-m)}H_{2}\left(
X_{m}\right)  $ for $m<n$ we get:
\begin{align*}
\mathbb{E}\left(  X_{n}^{2}X_{m}^{2}\right)   &  =\mathbb{E}X_{m}^{2}\left(
H_{2}\left(  X_{n}|q\right)  +1\right)  =1+\mathbb{E}\left(  X_{m}%
^{2}\mathbb{E}\left(  H_{2}\left(  X_{n}|q\right)  |\mathcal{F}_{\leq
m}\right)  \right) \\
&  =1+\rho^{2\left\vert n-m\right\vert }\mathbb{E}X_{m}^{2}H_{2}\left(
X_{m}\right)  =1+\rho^{2\left\vert n-m\right\vert }\mathbb{E}\left(  X_{m}%
^{4}-X_{m}^{2}\right) \\
&  =1+\rho^{2\left\vert n-m\right\vert }\left(  2+q-1\right)  .
\end{align*}
To get (\ref{2_1_1}), we have:%
\begin{align*}
\mathbb{E}X_{n}^{2}X_{n-j}X_{n+k}  &  =\mathbb{E}\left(  X_{n}^{2}%
X_{n-j}\mathbb{E}\left(  X_{n+k}|\mathcal{F}_{\leq n}\right)  \right)
=\rho^{k}\mathbb{E}\left(  X_{n}^{3}X_{n-j}\right) \\
&  =\rho^{k+j}\mathbb{E}\left(  X_{n}^{3}\mathbb{E}\left(  X_{n-j}%
|\mathcal{F}_{\geq n}\right)  \right)  =\rho^{k+j}\mathbb{E}X_{n}^{4}%
=\rho^{k+j}\left(  2+q\right)  .
\end{align*}

\end{proof}

\begin{proof}
[Proof of Proposition \ref{u_bryc2}]For fixed natural numbers $k,j$ let us
denote $Z_{i}\allowbreak=$\newline$\allowbreak\mathbb{E}\left(  X_{n+i}%
|\mathcal{F}_{\leq n-k,\geq n+j}\right)  $. Let $\mathbf{Z=}\left[
Z_{i}\right]  _{i=-k+1}^{j-1}$. Now notice that for each coordinate of the
vector $\mathbf{Z}$ we have:%
\begin{align*}
Z_{i}  &  =\mathbb{E}\left(  X_{i}|\mathcal{F}_{\leq n-k,\geq n+j}\right)
=\mathbb{E}\left(  \mathbb{E}\left(  X_{i}|\mathcal{F}_{\neq i}\right)
|\mathcal{F}_{\leq n-k,\geq n+j}\right) \\
&  =\frac{\rho}{1+\rho^{2}}\left(  \mathbb{E}\left(  X_{i-1}|\mathcal{F}_{\leq
n-k,\geq n+j}\right)  +\mathbb{E}\left(  X_{i+1}|\mathcal{F}_{\leq n-k,\geq
n+j}\right)  \right) \\
&  =\left\{
\begin{array}
[c]{ccc}%
\frac{\rho}{1+\rho^{2}}\left(  Z_{i-1}+Z_{i+1}\right)  & if & i\neq-k+1\vee
j-1\\
\frac{\rho}{1+\rho^{2}}X_{n-k}+\frac{\rho}{1+\rho^{2}}Z_{-k+2} & if & i=-k+1\\
\frac{\rho}{1+\rho^{2}}X_{n+j}+\frac{\rho}{1+\rho^{2}}Z_{j-2} & if & i=j-1
\end{array}
\right.  .
\end{align*}
Hence we have a vector linear equation
\[
\mathbf{Z=AZ+J},
\]
where
\begin{align*}
\mathbf{A}  &  \mathbf{=}\left[
\begin{array}
[c]{cccc}%
0 & \frac{\rho}{1+\rho^{2}} & \ldots & 0\\
\frac{\rho}{1+\rho^{2}} & 0 & \ldots & 0\\
\ldots & \ldots & \ldots & \frac{\rho}{1+\rho^{2}}\\
0 & 0 & \ldots & 0
\end{array}
\right]  ,\\
\mathbf{J}^{T}  &  =[\frac{\rho}{1+\rho^{2}}X_{n-k},0,\ldots,0,\frac{\rho
}{1+\rho^{2}}X_{n+j}].
\end{align*}
Now notice that the matrix $\mathbf{I-A}$ where $\mathbf{\ I}$ denotes unity
matrix is non-singular. This is so because sum of the absolute values of
elements in each row of the matrix $\mathbf{A}$ is less than $1$ which means
that eigenvalues of the matrix $\mathbf{I-A}$ are inside circle at center in
$1$ and radius less than $1,$ thus nonzero. Consequently each component of the
vector $\mathbf{\ Z}$ is a linear function of $X_{n-k}$ and $X_{n+j}$. Having
linearity of $\mathbb{E}\left(  X_{i}|\mathcal{F}_{\leq n-k,\geq n+j}\right)
$ with respect to $X_{n-k}$ and $X_{n+j}$ we get (\ref{liniowosc}).

Let us denote $m_{i,j}=\mathbb{E}\left(  X_{n+i}X_{n+j}|\mathcal{F}_{\leq
n-l,\geq n+k}\right)  ,$ for $i,j=-l,-l+1,\ldots,k-1,k$. Notice that using
(\ref{liniowosc}) we get $m_{i,j}=m_{j,i}$ and that $m_{-l,j}\allowbreak
=\allowbreak X_{n-l}\mathbb{E}\left(  X_{n+j}|\mathcal{F}_{\leq n-l,\geq
n+k}\right)  \allowbreak=\allowbreak\frac{\rho^{l+j}\left(  1-\rho
^{2k-2j}\right)  }{1-\rho^{2\left(  k+l\right)  }}X_{n-l}^{2}+\frac{\rho
^{k-j}\left(  1-\rho^{2l+2j}\right)  }{1-\rho^{2\left(  k+l\right)  }}%
X_{n-l}X_{n+k}$ and $m_{k,j}\allowbreak=\allowbreak X_{n+k}\mathbb{E}\left(
X_{n+j}|\mathcal{F}_{n-l\leq,\geq n+k}\right)  \allowbreak=\allowbreak
\frac{\rho^{l+j}\left(  1-\rho^{2k-2j}\right)  }{1-\rho^{2\left(  k+l\right)
}}X_{n-l}X_{n+k}\allowbreak+\allowbreak\frac{\rho^{k-j}\left(  1-\rho
^{2l+2j}\right)  }{1-\rho^{2\left(  k+l\right)  }}X_{n+k}^{2},$ $m_{-l,-l}%
=X_{n-1}^{2},$ $m_{k,k}=X_{n+k}^{2}$. Besides we have for $i,j=-j+1,\ldots
,k-1$ and $i\neq j$%
\begin{align*}
m_{i,j}  &  =\mathbb{E}\left(  X_{n+i}X_{n+j}|\mathcal{F}_{\leq n-l,\geq
n+k}\right)  =\mathbb{E}\left(  X_{n+i}\mathbb{E}\left(  X_{n+j}%
|\mathcal{F}_{\neq n+j}\right)  |\mathcal{F}_{n-l\leq,\geq n+k}\right) \\
&  =\frac{\rho}{1+\rho^{2}}\left(  m_{i,j-1}+m_{i,j+1}\right)
\end{align*}
and if $i=j$
\begin{align*}
m_{i,i}  &  =\mathbb{E}\left(  X_{n+i}^{2}|\mathcal{F}_{\leq n-l,\geq
n+k}\right)  =\mathbb{E}\left(  \mathbb{E}\left(  X_{n+i}^{2}|\mathcal{F}%
_{\neq n+i}\right)  |\mathcal{F}_{n-l\leq,\geq n+k}\right) \\
&  =A\left(  m_{i-1,i-1}+m_{i+1,i+1}\right)  +Bm_{i-1,i+1}+C
\end{align*}
Notice also that we have in fact $\left(  l+k-2\right)  ^{2}$ unknowns and
$\left(  l+k-2\right)  ^{2}$ linear equations. Moreover random variables
$m_{ij}$ are well defined since conditional expectation is uniquely defined
(up to set of probability $1)$. Thus we get the main assertion of the
proposition. Now we know that for some $A_{j},$ $B_{j},$ $C_{j}$ we have
\begin{equation}
\mathbb{E}\left(  X_{n}^{2}|\mathcal{F}_{n-j\leq,\geq n+k}\right)
=A_{jk}^{(1)}X_{n-j}^{2}+A_{jk}^{\left(  2\right)  }X_{n+k}^{2}+B_{jk}%
X_{n-j}X_{n+k}+C_{jk}. \label{rownosc}%
\end{equation}
First thing to notice is that (consequence of calculating expectation of both
sides of (\ref{rownosc}))
\[
1=A_{jk}^{\left(  1\right)  }+A_{jk}^{\left(  2\right)  }+\rho^{j+k}%
B_{jk}+C_{jk}.
\]
Secondly let us multiply (\ref{rownosc}) by $X_{n-j}^{2}$, $X_{n+k}^{2}$ and
$X_{n-j}X_{n+k}$ and let us take expectation of both sides of obtained in that
way equalities. In doing so we apply assertions of Lemma \ref{kwadraty}. In
this way we will get three equations:
\begin{align*}
1+\rho^{2j}(1+q)  &  =A_{jk}^{\left(  1\right)  }(2+q)+A_{jk}^{\left(
2\right)  }\left(  1+\rho^{2\left(  j+k\right)  }\left(  1+q\right)  \right)
\\
&  +B_{jk}\rho^{j+k}\left(  2+q\right)  +1-A_{jk}^{\left(  1\right)  }%
-A_{jk}^{\left(  2\right)  }-\rho^{j+k}B_{jk}\\
1+\rho^{2k}(1+q)  &  =A_{jk}^{\left(  1\right)  }\left(  1+\rho^{2\left(
j+k\right)  }\left(  1+q\right)  \right)  +A_{jk}^{\left(  2\right)  }(2+q)\\
&  +B_{jk}\rho^{j+k}\left(  2+q\right)  +1-A_{jk}^{\left(  1\right)  }%
-A_{jk}^{\left(  2\right)  }-\rho^{j+k}B_{jk}\\
\rho^{j+k}\left(  2+q\right)   &  =A_{jk}^{\left(  1\right)  }\rho
^{j+k}\left(  2+q\right)  +A_{jk}^{\left(  2\right)  }\rho^{j+k}\left(
2+q\right)  +B_{jk}\left(  1+\rho^{2\left(  j+k\right)  }\left(  1+q\right)
\right) \\
&  +\rho^{j+k}\left(  1-A_{jk}^{\left(  1\right)  }-A_{jk}^{\left(  2\right)
}-\rho^{j+k}B_{jk}\right)  .
\end{align*}
Solution of this system of equations is (\ref{par_j}), as it can be easily
checked .
\end{proof}

\begin{proof}
[Proof of lemma \ref{general}]Consider discrete time random field
$\mathbf{Z=}\left\{  Z_{k}\right\}  _{k\in\mathbb{Z}}$ such that
$Z_{k}=X_{kj+m},$ for some fixed $j$ and $0\leq m\leq j-1$. Obviously we have
\begin{align*}
\mathbb{E}\left(  Z_{k}|\mathcal{F}_{\neq k}\right)   &  =\mathbb{E}\left(
X_{kj+m}|\mathcal{F}_{\leq\left(  k-1\right)  j+m,\geq\left(  k+1\right)
j+m}\right) \\
&  =\frac{\rho^{j}}{1+\rho^{2j}}\left(  X_{\left(  k-1\right)  j+m}+X_{\left(
k+1\right)  j+m}\right)  =\frac{\rho^{j}}{1+\rho^{2j}}\left(  Z_{k-1}%
+Z_{k+1}\right)  ,
\end{align*}
and
\begin{align*}
\mathbb{E}\left(  Z_{k}^{2}|\mathcal{F}_{\neq k}\right)   &  =\mathbb{E}%
\left(  X_{kj+m}^{2}|\mathcal{F}_{\leq\left(  k-1\right)  j+m,\geq\left(
k+1\right)  j+m}\right) \\
&  =A_{j}\left(  X_{\left(  k-1\right)  j+m}^{2}+X_{\left(  k+1\right)
j+m}^{2}\right)  +B_{j}X_{\left(  k-1\right)  j+m}X_{\left(  k+1\right)
j+m}+C_{j}\\
&  =A_{j}\left(  Z_{k-1}^{2}+Z_{k+1}^{2}\right)  +B_{j}Z_{k-1}Z_{k+1}+C_{j},
\end{align*}
where $A_{j}=A_{jj}^{\left(  1\right)  }=A_{jj}^{\left(  2\right)  },$
$B_{j}=B_{jj}$. Thus $\mathbf{Z}$ is 1TSP with different parameters. Notice
that one dimensional distributions of processes $\mathbf{X}$ and $\mathbf{Z}$
are the same. Hence parameters $q$ for both processes $\mathbf{Z}$ and
$\mathbf{X}$ are the same. On the other hand parameter $\rho_{Z}$ of the
process $\mathbf{Z}$ is related to parameter $\rho$ of the process
$\mathbf{X}$ by the following relationship
\[
\rho_{Z}=\mathbb{E}Z_{0}Z_{1}=\mathbb{E}X_{0}X_{j}=\rho^{j}.
\]
Thus applying formulae (\ref{_B}) we get (\ref{par_j}).
\end{proof}

\begin{proof}
[Proof of Theorem \ref{existence}]From Lemma \ref{general} of the previous
section it follows that if for some $\delta>0$ $\mathbf{X\allowbreak
=\allowbreak}\left\{  X_{n}^{\left(  \delta\right)  }\right\}  _{n\in
\mathbb{Z}}$ is a regular 1TSP with some parameters $q$ and $\rho,$ then for
every $j$ and $m$\allowbreak$\in$\allowbreak$\left\{  0,\ldots,j-1\right\}  ,$
process $Z_{k}^{\left(  m\right)  }\allowbreak=\allowbreak X_{kj+m}^{\left(
\delta\right)  }$ is also the 1TSP with the same parameter $q$ and parameter
$\rho_{z}\allowbreak=\allowbreak\rho^{j}.$ Notice however that $Z_{k}%
^{(0)}\allowbreak=\allowbreak X_{k}^{\left(  j\delta\right)  }.$ Similarly if
we considered process $\mathbf{\hat{X}\allowbreak=\allowbreak}\left\{  \hat
{X}_{n}^{\left(  \delta/j\right)  }\right\}  _{n\in\mathbb{Z}},$ then since
$X_{n}^{\left(  \delta\right)  }\allowbreak=\allowbreak\hat{X}_{nj}^{\left(
\delta/j\right)  }$ and the fact that one 1TSP is characterized by one
parameter $q$ we deduce that processes $\mathbf{X}$ and $\mathbf{\hat{X}}$
share the same parameter $q.$ Hence the fact that condition $\left(
\delta\right)  $ is satisfied for every $\delta$ implies that all implied by
it regular 1TSP are characterized by one, same parameter $q.$ Further since
for same $\delta>0\mathbf{\ }$regular 1TSP $\left\{  X_{n}^{\left(
\delta\right)  }\right\}  _{n\in\mathbb{Z}}$ is $L_{2}$ stationary with
covariance function $K\left(  n,m\right)  \allowbreak=\allowbreak\rho\left(
\delta\right)  ^{\left\vert n-m\right\vert }\allowbreak\mathbb{=\allowbreak
\mathbb{E}}X_{0}^{\left(  \delta\right)  }X_{\left\vert n-m\right\vert
}^{\left(  \delta\right)  }$, we have also for any integer $k$: $\rho\left(
\delta\right)  ^{n}\mathbb{\allowbreak=\allowbreak\mathbb{E}}X_{0}^{\left(
\delta\right)  }X_{n}^{\left(  \delta\right)  }\allowbreak=\allowbreak
\mathbb{E}X_{0}^{\left(  \delta/k\right)  }X_{nk}^{\left(  \delta/k\right)
}\allowbreak=\allowbreak\rho\left(  \delta/k\right)  ^{nk}.$ Or equivalently
we have
\[
\forall k\in\mathbb{Z};\delta,\theta>0:\rho\left(  \delta\right)
\allowbreak=\allowbreak\rho\left(  \delta/k\right)  ^{k},\rho\left(
k\theta\right)  \allowbreak=\allowbreak\rho\left(  \theta\right)  ^{k}%
\]
Now take $\delta\allowbreak=\allowbreak\frac{k}{m}\theta$ for some $\theta.$
We will get then $\rho\left(  \frac{k}{m}\theta\right)  \allowbreak
=\allowbreak\rho\left(  \theta\right)  ^{\frac{k}{m}}\allowbreak.$ Now let us
take sequences of integers $\left\{  k_{n},m_{n}\right\}  $ such that
$\frac{k_{n}}{m_{n}}\underset{n\rightarrow\infty}{\longrightarrow}1/\theta$
then, using $L_{2}-$ continuity we get%
\[
\rho\left(  \theta\right)  \allowbreak=\allowbreak\rho\left(  1\right)
^{\theta}.
\]
In other words we deduce that if $\mathbf{Y}$ existed, then it would be
$L_{2}-$ stationary with covariance function
\[
K\left(  t,s\right)  =K\left(  \left\vert s-t\right\vert \right)  =\rho\left(
1\right)  ^{\left\vert s-t\right\vert }.
\]
for some $\rho\left(  1\right)  \in(0,1)$. Let us introduce new parameter
\[
\alpha=\log\frac{1}{\rho\left(  1\right)  }>0.
\]
We have then
\[
K\left(  s,t\right)  =\exp\left(  -\alpha\left\vert s-t\right\vert \right)  .
\]
Consequently the fact that condition $\left(  \delta\right)  $ is satisfied
for every $\delta$ implies that all implied by it regular 1TSP are
characterized by one, same parameter $q$ and covariance function defined by
same parameter $\alpha.$

Existence of $\mathbf{Y}$ will be shown for two cases separately. Since for
$q\allowbreak=\allowbreak1$ we have normality of the one dimensional and
conditional distributions. Thus we deduce that the process $\mathbf{Y}$ for
$q\allowbreak=\allowbreak1$ is in fact the well known Ornstein--Uhlenbeck process.

Now let us consider fixed $q\in(-1,1)$. First we will deduce the existence of
the process $\mathbf{\tilde{Y}=}\left(  Y_{t}\right)  _{t\in\mathbb{Q}}$. This
follows from Kolmogorov's extension theorem. Since having natural ordering of
$\mathbb{Q}$ we need only consistency of the family of finite dimensional
distributions of $\mathbf{\tilde{Y}}$. This can be however easily shown by the
following argument. Let us take finite set of numbers $r_{1}<r_{2}%
<\cdots<r_{n}$ from $\mathbb{Q}$. Let $M$ denote the smallest common
denominator of these numbers. Let us consider regular 1TSP $\mathbf{X}^{M}$
with $q$ and $\hat{\alpha}=\alpha/M$. Note that then $X_{n}^{M}=Y_{n/M}$ for
$n\in\mathbb{N}$. Then joint distribution of $\left(  Y_{r_{1}},\ldots
,Y_{r_{n}}\right)  $ is in fact a joint distribution of $\left(  X_{R_{1}}%
^{M},\ldots,X_{R_{n}}^{M}\right)  $ where numbers $R_{1},\ldots,R_{n}$ are
defined by the relationships $r_{i}=R_{i}/M$. Since process $\mathbf{X}^{M}$
exists we have consistency since if $\left\{  \tau_{1},\ldots,\tau
_{k}\right\}  \subset\left\{  r_{1},\ldots,r\,_{n}\right\}  $ then
distribution of $\left(  Y_{\tau_{1}},\ldots,Y_{\tau_{k}}\right)  $ being
equal to the distribution of $\left(  X_{T_{1}}^{M},\ldots,X_{T_{k}}%
^{M}\right)  $ with $T_{j}$ defined by $\tau_{j}=T_{j}/M$ for $j=1,\ldots,k$
is a projection of the distribution of $\left(  Y_{r_{1}},\ldots,Y_{r_{n}%
}\right)  $. Hence we deduce that the process $\mathbf{\tilde{Y}}$ with values
in the compact space $S\left(  q\right)  ,$ exist. Now we use separability
theorem (see e.g. \cite{wencel}) and view $\mathbf{\tilde{Y}}$ as separable
modification of the process $\mathbf{Y}$ itself. Hence the process
$\mathbf{Y}$ exists.
\end{proof}

\begin{proof}
[Proof of Theorem \ref{ou-_wlasn}]i) The fact that $\mathbb{E}Y_{t}%
^{2}\allowbreak=\allowbreak1$ and that for $s\geq0,$ $\operatorname*{var}%
\left(  Y_{t+s}|Y_{t}\right)  \allowbreak=\allowbreak1-e^{-2\alpha s}$ follows
from (\ref{_hermi}) and (\ref{_alsalam}) and the definition of polynomials
$p_{n}$ for $n=1,2$. If $q\allowbreak=\allowbreak1$ than we have OU process
and the assertion is true. For $\left\vert q\right\vert <1,$ we apply
assertion of the Theorem 3 page 180 of \cite{wencel} and the following
estimation based on Chebyshev inequality.
\begin{gather*}
\gamma_{\varepsilon}\left(  t\right)  \allowbreak=\allowbreak\sup_{y\in
S\left(  q\right)  ,t\leq h}P\left(  |Y_{s+t}-Y_{s}|\geq\varepsilon
|Y_{s}=y\right)  \allowbreak\leq\\
\frac{\mathbb{E}\left\vert Y_{t+s}-Y_{s}\right\vert ^{2}}{\varepsilon^{2}%
}\allowbreak\leq\allowbreak\frac{2\operatorname*{var}\left(  Y_{t+s}%
|Y_{s}=y\right)  +2\mathbb{E(}\left\vert Y_{s}-e^{-\alpha t}Y_{s}\right\vert
^{2}|Y_{s}=y)}{\varepsilon^{2}}\\
=\allowbreak\frac{2\left(  1-e^{-2\alpha t}\right)  +2y^{2}\left(
1-e^{-\alpha t}\right)  ^{2}}{\varepsilon^{2}}\allowbreak\cong\allowbreak
\frac{4\alpha t+8\alpha^{2}t^{2}/(1-q)}{\varepsilon^{2}}%
\end{gather*}
$\allowbreak\allowbreak$ from which it follows that $\forall\varepsilon>0,$
$\gamma_{\varepsilon}\left(  t\right)  ->0$ as $t->0$. Another justification
of this assertion follows properties of martingales and is given below.

ii) Following observations: 1. $q-$Hermite polynomials are the orthogonal
basis of the space denoted by $B\left(  q\right)  $. 2. Conditional
distributions of $Y_{t+s}|Y_{t}=y$ having densities $f_{CN}\left(
x|y,e^{-\alpha s},q\right)  $ form a continuous semigroup following Lemma
\ref{Hip}, iv). 3. (\ref{formula}) follows Lemma \ref{Hip}, ii).

iii) If $q\allowbreak=\allowbreak1$ then we deal with classical OU process
that has both Feller property and is strongly Markov. For $\left\vert
q\right\vert <1$ we use the fact that $\left\vert H_{i}\left(  x|q\right)
\right\vert \allowbreak\leq\allowbreak W_{i}\left(  q\right)  /(1-q)^{i/2}$ by
Lemma \ref{Hip}, v), (\ref{formula}) together with Lemma \ref{Hip}, vi)
guarantees that
\[
\max_{x\in S\left(  q\right)  }\left\vert P^{\tau}(f)(x)\right\vert \leq
\frac{\max_{j\geq0}\left\vert b_{j}\sqrt{\left(  q\right)  _{j}}\right\vert
}{\left(  \exp\left(  -\alpha\tau\right)  \right)  _{\infty}^{2}},
\]
where coefficients $b_{i}$ are defined by
\[
B\left(  q\right)  \ni f\left(  x\right)  \allowbreak=\allowbreak\sum_{j\geq
0}\frac{b_{j}}{\sqrt{\left[  j\right]  _{q}!}}H_{j}\left(  x|q\right)  .
\]
Thus if $f$ was continuous and bounded, then $P^{\tau}\left(  f\right)
\left(  x\right)  $ is also continuous and bounded. Hence we have Feller
property. To get strong Markov property we use Theorem $1$ of section 9.2 of
\cite{wencel} that asserts that every time homogeneous Markov family, having
c\'{a}dl\'{a}g trajectories and Feller property is also strongly Markov.

iv) Let us consider function $f\allowbreak\in\allowbreak B\left(  q\right)  $
and take $n_{0}\allowbreak\in\allowbreak\mathbb{N}$ such that%
\[
\max(\sup_{x\in S\left(  q\right)  }\left\vert \sum_{n\geq n_{0}}\frac{b_{n}%
}{\sqrt{\left[  n\right]  _{q}!}}H_{n}\left(  x|q\right)  \right\vert
,\sup_{x\in S\left(  q\right)  }\left\vert \sum_{n\geq n_{0}}e^{-n\alpha\tau
}\frac{b_{n}}{\sqrt{\left[  n\right]  _{q}!}}H_{n}\left(  x|q\right)
\right\vert \allowbreak)\leq\varepsilon
\]
for some chosen beforehand $\varepsilon>0.$ Since on the compact space
$S\left(  q\right)  $ uniform convergence to a continuous function is
equivalent to pointwise convergence and since from $L_{2}$ convergence follows
existence of a subsequence $\left\{  k_{n}\right\}  $ such that $\sum_{j\geq
0}^{k_{n}}\frac{b_{j}}{\sqrt{\left[  j\right]  _{q}!}}H_{j}\left(  x|q\right)
$ converges pointwise to its continuous limit, we deduce that such $n_{0}$
exists. Now we notice
\[
\sup_{x\in\mathbb{S}\left(  q\right)  }\left\vert f\left(  x\right)  -P^{\tau
}\left(  f\right)  \left(  x\right)  \right\vert \leq2\varepsilon+\sup_{x\in
S\left(  q\right)  }\left\vert \sum_{j=0}^{n_{0}-1}\left(  1-e^{-j\alpha\tau
}\right)  \frac{b_{j}}{\sqrt{\left[  j\right]  _{q}!}}H_{j}\left(  x|q\right)
\right\vert \leq3\varepsilon,
\]
if only $\tau$ is sufficiently small. Hence
\[
\lim_{\tau\rightarrow0^{+}}\sup_{x\in\mathbb{S}\left(  q\right)  }\left\vert
f\left(  x\right)  -P^{\tau}\left(  f\right)  \left(  x\right)  \right\vert
=0,
\]
and consequently we see that $\left\{  P^{\tau}\left(  f\right)  \right\}
_{\tau\geq0}$ is a right continuous meaning that process $\mathbf{Y}$ is a
Feller process. Thus the infinitesimal operator $A$ exists. Its value on
$H_{n}\left(  x|q\right)  $ can be found by the following argument:
\begin{align*}
AH_{n}\left(  x|q\right)  \allowbreak &  =\allowbreak\lim_{h\downarrow0^{+}%
}\frac{1}{h}(\mathbb{E}\left(  H_{n}\left(  Y_{t+h}|q\right)  |Y_{t}=x\right)
\allowbreak-\allowbreak H_{n}\left(  x|q)\right)  \allowbreak\\
&  =\allowbreak\lim_{h\downarrow0^{+}}\frac{1}{h}\left(  e^{-\alpha nh}%
H_{n}\left(  x|q\right)  -H_{n}\left(  x|q\right)  \right)  \allowbreak
=\allowbreak-n\alpha H_{n}\left(  x|q\right)  .
\end{align*}

\end{proof}

\begin{proof}
[Proof of Theorem \ref{wl_q-Wienera}]1. Suppose $\tau<\sigma,$
\begin{align*}
\operatorname{cov}\left(  X_{\tau},X_{\sigma}\right)  \allowbreak &
=\allowbreak\sqrt{\tau\sigma}\mathbb{E}Y_{\log\tau/2\alpha}Y_{\log
\sigma/2\alpha}\allowbreak=\allowbreak\sqrt{\tau\sigma}e^{-\alpha\left\vert
\log\tau/2\alpha-\log\sigma/2\alpha\right\vert }\allowbreak\\
&  =\allowbreak\sqrt{\tau\sigma}e^{-\log\left(  \sigma/\tau\right)
/2}\allowbreak=\allowbreak\sqrt{\tau\sigma}\sqrt{\frac{\sigma}{\tau}}=\sigma.
\end{align*}
2. \& 3. We have $\mathbb{P}\left(  X_{t}\leq x\right)  =\allowbreak
\mathbb{P}\left(  \sqrt{\tau}Y_{\log\tau/2\alpha}\leq x\right)  \allowbreak
=\allowbreak\mathbb{P}\left(  Y_{\log\tau/2\alpha}\leq x/\sqrt{\tau}\right)
$. Knowing that $f_{N}$ is a density of $Y_{t}$ we get immediately assertion
1. To get assertion 3 we have
\begin{align*}
\mathbb{P}\left(  X_{\tau}-X_{\sigma}\leq x|X_{\sigma}=y\right)  \allowbreak
&  =\allowbreak\mathbb{P}\left(  X_{\tau}\leq x+y|X_{\sigma}=y\right)
\allowbreak=\\
\mathbb{P}\left(  \sqrt{\tau}Y_{\log\tau/2\alpha}\leq x+y|\sqrt{\sigma}%
Y_{\log\sigma/2\alpha}=y\right)  \allowbreak &  =\allowbreak\mathbb{P}\left(
Y_{\log\tau/2\alpha}\leq\frac{x+y}{\sqrt{\tau}}|Y_{\log\sigma/2\alpha}%
=\frac{y}{\sqrt{\sigma}}\right)  .
\end{align*}
Now we recall that $f_{CN}\left(  x|y,q,e^{-\alpha\left\vert s-t\right\vert
}\right)  $ is the density of $Y_{t}$ given $Y_{s}=s$, by theorem
\ref{podstawowe}.

4. Notice that $X_{\tau}/\sqrt{\tau}\allowbreak=\allowbreak Y_{\log
\tau/2\alpha}$. We have using assertion $5$ of Theorem \ref{podstawowe}:
\begin{align*}
\mathbb{E}\left(  \tau^{n/2}H_{n}\left(  X_{\tau}/\sqrt{\tau}\right)
|\mathcal{F}_{\leq\sigma}^{X}\right)  \allowbreak &  =\allowbreak\tau
^{n/2}\mathbb{E}\left(  H_{n}\left(  Y_{\log\tau/2\alpha}\right)
|\mathcal{F}_{\log\sigma/2\alpha}\right)  \allowbreak=\allowbreak\\
\tau^{n/2}e^{-\alpha n\left(  \log\tau/2\alpha-\log\sigma/2\alpha\right)
}H_{n}\left(  Y_{\log\sigma/2\alpha}\right)  \allowbreak &  =\allowbreak
\sigma^{n/2}H_{n}\left(  X_{\sigma}/\sqrt{\sigma}\right)  .
\end{align*}
a.s. \newline Similarly
\begin{align*}
\mathbb{E}\left(  \sigma^{-n/2}H_{n}\left(  \frac{X_{\sigma}}{\sqrt{\sigma}%
}|q\right)  |\mathcal{F}_{\geq\tau}^{X}\right)  \allowbreak &  =\allowbreak
\sigma^{-n/2}\mathbb{E}\left(  H_{n}\left(  Y_{\log\sigma/2\alpha}|q\right)
|\mathcal{F}_{\geq\log\tau/2\alpha}\right)  \allowbreak=\\
\sigma^{-n/2}e^{-\alpha\left(  \log\tau/2\alpha-\log\sigma/2\alpha\right)
}H_{n}\left(  Y_{\log\tau/2\alpha}|q\right)  \allowbreak &  =\newline%
\allowbreak\tau^{-n/2}H_{n}\left(  \frac{X_{\tau}}{\sqrt{\tau}}|q\right)  .
\end{align*}
\newline5. to get first part we proceed as follows:
\begin{align*}
\sqrt{\sigma}\mathbb{E}\left(  Y_{\log\sigma/2\alpha}|\mathcal{F}_{\leq
\log\left(  \sigma-\delta\right)  /2\alpha,\geq\log\left(  \sigma
+\gamma\right)  /2\alpha}\right)  \allowbreak &  =\sqrt{\sigma}(\frac
{\sqrt{\frac{\left(  \sigma-\delta\right)  }{\sigma}}\left(  1-\sigma/\left(
\sigma+\gamma\right)  \right)  }{1-\frac{\left(  \sigma-\delta\right)
}{\left(  \sigma+\gamma\right)  }}Y_{\log\left(  \sigma-\delta\right)
/2\alpha}\allowbreak\\
+\allowbreak\frac{\sqrt{\frac{\sigma}{\left(  \sigma+\gamma\right)  }}\left(
1-\left(  \sigma-\delta\right)  /\sigma\right)  }{1-\frac{\left(
\sigma-\delta\right)  }{\left(  \sigma+\gamma\right)  }}Y_{\log\left(
\sigma+\gamma\right)  /2\alpha})  &  =\allowbreak\frac{\gamma}{\delta+\gamma
}X_{\sigma-\delta}+\frac{\delta}{\delta+\gamma}X_{\sigma+\gamma}.
\end{align*}
$\allowbreak$ To get second part we have:
\begin{gather*}
\mathbb{E}\left(  X_{\sigma}^{2}|\mathcal{F}_{\leq\sigma-\delta,\geq
\sigma+\gamma}^{X}\right)  \allowbreak=\allowbreak\newline\sigma
\mathbb{E}\left(  Y_{\log\sigma/2\alpha}^{2}|\mathcal{F}_{\leq\log\left(
\sigma-\delta\right)  /2\alpha,\geq\log\left(  \sigma+\gamma\right)  /2\alpha
}\right)  \allowbreak\\
=\sigma\frac{\frac{\sigma-\delta}{\sigma}\left(  1-\frac{\sigma}{\sigma
+\gamma}\right)  \left(  1-\frac{q\sigma}{\sigma+\gamma}\right)  }{\left(
1-\frac{\sigma-\delta}{\sigma+\gamma}\right)  \left(  1-q\frac{\sigma-\delta
}{\sigma+\gamma}\right)  }Y_{\log\left(  \sigma-\delta\right)  /2\alpha}%
^{2}\allowbreak+\allowbreak\sigma\frac{\frac{\sigma}{\sigma+\gamma}\left(
1-\frac{\sigma-\delta}{\sigma}\right)  \left(  1-q\frac{\sigma-\delta}{\sigma
}\right)  }{\left(  1-\frac{\sigma-\delta}{\sigma+\gamma}\right)  \left(
1-q\frac{\sigma-\delta}{\sigma+\gamma}\right)  }Y_{\log\left(  \sigma
+\gamma\right)  /2\alpha}^{2}\\
+\allowbreak\sigma\frac{\left(  q+1\right)  \sqrt{\frac{\sigma-\delta}%
{\sigma+\gamma}}\left(  1-\frac{\sigma}{\sigma+\gamma}\right)  \left(
1-\frac{\sigma-\delta}{\sigma}\right)  }{\left(  1-\frac{\sigma-\delta}%
{\sigma+\gamma}\right)  \left(  1-q\frac{\sigma-\delta}{\sigma+\gamma}\right)
}Y_{\log\left(  \sigma-\delta\right)  /2\alpha}Y_{\log\left(  \sigma
+\gamma\right)  /2\alpha}\allowbreak+\allowbreak\allowbreak\sigma\frac{\left(
1-\frac{\sigma-\delta}{\sigma}\right)  \left(  1-\frac{\sigma}{\sigma+\gamma
}\right)  }{\left(  1-q\frac{\sigma-\delta}{\sigma+\gamma}\right)  }.
\end{gather*}

$\allowbreak$Keeping in mind that $X_{\tau}=\sqrt{\tau}Y_{\log\tau/2\alpha}$
we get:
\begin{gather*}
\mathbb{E}\left(  X_{\sigma}^{2}|\mathcal{F}_{\leq\sigma-\delta,\geq
\sigma+\gamma}^{X}\right)  \allowbreak=\allowbreak\newline\frac{\gamma\left(
\left(  1-q\right)  \sigma+\gamma\right)  }{\left(  \delta+\gamma\right)
\left(  \sigma\left(  1-q\right)  +\gamma+q\delta\right)  }X_{\sigma-\delta
}^{2}\allowbreak+\allowbreak\frac{\delta\left(  \left(  1-q\right)
\sigma+q\delta\right)  }{\left(  \delta+\gamma\right)  \left(  \sigma\left(
1-q\right)  +\gamma+q\delta\right)  }X_{\sigma+\gamma}^{2}\allowbreak\\
+\allowbreak\frac{\left(  q+1\right)  \delta\gamma}{\left(  \delta
+\gamma\right)  \left(  \sigma\left(  1-q\right)  +\gamma+q\delta\right)
}X_{\sigma-\delta}X_{\sigma+\gamma}\allowbreak+\allowbreak\frac{\delta\gamma
}{\left(  \sigma\left(  1-q\right)  +\gamma+q\delta\right)  }.
\end{gather*}

\end{proof}

\begin{proof}
[Proof of Corollary \ref{har4}]To prove $i)$ we calculate:%

\begin{align*}
\mathbb{E}\left(  \left(  X_{\tau}-X_{\sigma}\right)  ^{2}|\mathcal{F}%
_{\leq\sigma}^{X}\right)   &  =\mathbb{E}\left(  X_{\tau}^{2}|\mathcal{F}%
_{\leq\sigma}^{X}\right)  -2X_{\sigma}^{2}+X_{\sigma}^{2}=X_{\sigma}^{2}%
+\tau-\sigma-X_{\sigma}^{2}=\tau-\sigma.\\
\mathbb{E}\left(  \left(  X_{\tau}-X_{\sigma}\right)  ^{3}|\mathcal{F}%
_{\leq\sigma}^{X}\right)  \allowbreak &  =\allowbreak\mathbb{E}\left(
X_{\tau}^{3}|\mathcal{F}_{\leq\sigma}^{X}\right)  -3X_{\sigma}\mathbb{E}%
\left(  X_{\tau}^{2}|\mathcal{F}_{\leq\sigma}^{X}\right)  \allowbreak
+\allowbreak3X_{\sigma}^{2}\mathbb{E}\left(  X_{\tau}|\mathcal{F}_{\leq\sigma
}^{X}\right)  -X_{\sigma}^{3}\allowbreak\\
&  =\allowbreak\mathbb{E}\left(  X_{\tau}^{3}|\mathcal{F}_{\leq\sigma}%
^{X}\right)  -3X_{\sigma}\mathbb{E}\left(  X_{\tau}^{2}|\mathcal{F}%
_{\leq\sigma}^{X}\right)  \allowbreak+\allowbreak2X_{\sigma}^{3}.
\end{align*}

\begin{align*}
\mathbb{E}\left(  \left(  X_{\tau}-X_{\sigma}\right)  ^{4}|\mathcal{F}%
_{\leq\sigma}^{X}\right)   &  =\mathbb{E}\left(  X_{\tau}^{4}|\mathcal{F}%
_{\leq\sigma}^{X}\right)  -4X_{\sigma}\mathbb{E}\left(  X_{\tau}%
^{3}|\mathcal{F}_{\leq\sigma}^{X}\right)  +6X_{\sigma}^{2}\mathbb{E}\left(
X_{\tau}^{2}|\mathcal{F}_{\leq\sigma}^{X}\right)  -4X_{\sigma}^{3}%
\mathbb{E}\left(  X_{\tau}|\mathcal{F}_{\leq\sigma}^{X}\right)  +X_{\sigma
}^{4}\\
&  =\mathbb{E}\left(  X_{\tau}^{4}|\mathcal{F}_{\leq\sigma}^{X}\right)
-4X_{\sigma}\mathbb{E}\left(  X_{\tau}^{3}|\mathcal{F}_{\leq\sigma}%
^{X}\right)  +3X_{\sigma}^{4}+6\left(  \tau-\sigma\right)  X_{\sigma}^{2}.
\end{align*}

Now recalling that $H_{2}\left(  x|q\right)  \allowbreak\allowbreak
=\allowbreak x^{2}-1,$ $\allowbreak H_{3}\left(  x|q\right)  \allowbreak
=\allowbreak x^{3}-(2+q)x$ and $H_{4}\left(  x|q\right)  \allowbreak
=\allowbreak x^{4}\allowbreak-\allowbreak(3+2q+q^{2})x^{2}\allowbreak
+\allowbreak(1+q+q^{2})$ we see that:
\begin{align*}
\sigma H_{2}\left(  \frac{X_{\sigma}}{\sqrt{\sigma}}|q\right)  \allowbreak &
=\allowbreak\mathbb{E}\left(  \tau H_{2}\left(  \frac{X_{\tau}}{\sqrt{\tau}%
}|q\right)  |\mathcal{F}_{\leq\sigma}^{X}\right)  ,\\
\allowbreak\sigma^{3/2}H_{3}\left(  \frac{X_{\sigma}}{\sqrt{\sigma}}|q\right)
\allowbreak &  =\allowbreak\mathbb{E}\left(  \tau^{3/2}H_{3}\left(
\frac{X_{\tau}}{\sqrt{\tau}}|q\right)  |\mathcal{F}_{\leq\sigma}^{X}\right) \\
&  =\allowbreak\mathbb{E}\left(  X_{\tau}^{3}|\mathcal{F}_{\leq\sigma}%
^{X}\right)  \allowbreak-\allowbreak\tau(2+q)\mathbb{E}\left(  X_{\tau
}|\mathcal{F}_{\leq\sigma}^{X}\right)
\end{align*}
$\allowbreak$ and%

\begin{align*}
\sigma^{2}H_{4}\left(  \frac{X_{\sigma}}{\sqrt{\sigma}}|q\right)  \allowbreak
&  =\allowbreak\mathbb{E}\left(  \tau^{2}H_{4}\left(  \frac{X_{t}}{\sqrt{\tau
}}|q\right)  |\mathcal{F}_{\leq\sigma}^{X}\right)  \allowbreak\\
&  =\allowbreak\mathbb{E}\left(  X_{\tau}^{4}|\mathcal{F}_{\leq\sigma}%
^{X}\right)  \allowbreak-\allowbreak\allowbreak\tau\left(  3+2q+q^{2}\right)
\allowbreak\times\allowbreak\mathbb{E}\left(  X_{\tau}^{2}|\mathcal{F}%
_{\leq\sigma}^{X}\right)  \allowbreak+\allowbreak(1+q+q^{2})\tau^{2}.
\end{align*}
Thus $\mathbb{E}\left(  X_{\tau}^{2}|\mathcal{F}_{\leq\sigma}^{X}\right)
\allowbreak=\allowbreak X_{\sigma}^{2}+\tau-\sigma$

$\allowbreak$%
\[
\mathbb{E}\left(  X_{\tau}^{3}|\mathcal{F}_{\leq\sigma}^{X}\right)
=\tau(2+q)X_{\sigma}+X_{\sigma}^{3}-\sigma\left(  2+q\right)  X_{\sigma
}\allowbreak=\allowbreak\left(  \tau-\sigma\right)  (2+q)X_{\sigma}+X_{\sigma
}^{3}%
\]
and
\[
\mathbb{E}\left(  X_{\tau}^{4}|\mathcal{F}_{\leq\sigma}^{X}\right)
\allowbreak=\allowbreak X_{\sigma}^{4}\allowbreak+\allowbreak\left(
\tau-\sigma\right)  (3+2q+q^{2})X_{\sigma}^{2}-\left(  \tau-\sigma\right)
\left(  \tau+\sigma\right)  (1+q+q^{2}).
\]

So%
\begin{align*}
\mathbb{E}\left(  \left(  X_{\tau}-X_{\sigma}\right)  ^{3}|\mathcal{F}%
_{\leq\sigma}^{X}\right)  \allowbreak &  =\allowbreak\left(  \tau
-\sigma\right)  (2+q)X_{\sigma}+X_{\sigma}^{3}-3X_{\sigma}\left(  X_{\sigma
}^{2}+\tau-\sigma\right)  +2X_{\sigma}^{3}\allowbreak\\
&  =\allowbreak-\left(  1-q\right)  \left(  \tau-\sigma\right)  X_{\sigma},
\end{align*}
\newline%
\begin{gather*}
\mathbb{E}\left(  \left(  X_{\tau}-X_{\sigma}\right)  ^{4}|\mathcal{F}%
_{\leq\sigma}^{X}\right)  \allowbreak=\allowbreak3X_{\sigma}^{4}+6\left(
\tau-\sigma\right)  X_{\sigma}^{2}\allowbreak+\allowbreak X_{\sigma}%
^{4}\allowbreak+\allowbreak\left(  \tau-\sigma\right)  \left(  3+2q+q^{2}%
\right)  X_{\sigma}^{2}\allowbreak\\
+\allowbreak\left(  \tau-\sigma\right)  \tau\left(  3+2q+q^{2}\right)
\allowbreak-\allowbreak(1+q+q^{2})\left(  \tau^{2}-\sigma^{2}\right)
\allowbreak-\allowbreak4X_{\sigma}^{4}-4\left(  \tau-\sigma\right)  \left(
2+q\right)  X_{\sigma}^{2}\\
\allowbreak=\allowbreak X_{\sigma}^{2}\left(  \tau-\sigma\right)  \left(
1-q\right)  ^{2}+\left(  2+q\right)  \left(  \tau-\sigma\right)  ^{2}+\left(
\tau-\sigma\right)  \sigma\left(  1-q^{2}\right)  .
\end{gather*}

If $\mathbf{X}$ had independent increments then $\mathbb{E}\left(  \left(
X_{\tau}-X_{\sigma}\right)  ^{4}|\mathcal{F}_{\leq\sigma}^{X}\right)  $ would
not depend on $X_{\sigma}^{2}.$ Note that for $q=1$ we have $\mathbb{E}\left(
\left(  X_{\tau}-X_{\sigma}\right)  ^{4}|\mathcal{F}_{\leq\sigma}^{X}\right)
\allowbreak=\allowbreak3\left(  \tau-\sigma\right)  ^{2}\allowbreak
=\allowbreak3\mathbb{E}\left(  \left(  X_{\tau}-X_{\sigma}\right)
^{2}|\mathcal{F}_{\leq\sigma}^{X}\right)  .$

ii) Follows properties of martingales. Thus almost every path of the process
$\mathbf{\ X}$ has no oscillatory discontinuities, in other words has at every
point $t$ left and right hand side limit. Since process $\mathbf{Y}$ is
obtained from the process $\mathbf{X}$ by continuous transformation it has
similar properties.

iii) Follows the fact that $q-$Wiener process is obtained from $(\alpha,q)-$OU
process by continuous (even smooth) transformation (\ref{def_q-W}). On its
side $(\alpha,q)-$OU has Feller property and is strongly Markov as shown in
Theorem \ref{ou-_wlasn},iii).
\end{proof}

\begin{proof}
[Proof of Corollary \ref{warobustr}]Following given formulae we have for $n=1$
and $n=2:$ $A_{0,0}^{(1)}\allowbreak=\allowbreak\frac{\gamma}{\gamma+\delta
}\sqrt{\frac{\sigma-\delta}{\sigma}},$ $A_{0,1}^{(1)}\allowbreak
=\allowbreak\frac{\delta}{\gamma+\delta}\sqrt{\frac{\sigma+\gamma}{\sigma}},$
$A_{0,-1}^{(2)}\allowbreak=\allowbreak\frac{\gamma\delta}{\sigma\left(
\gamma+\delta\right)  (\gamma+q\delta+(1-q)\sigma)}\frac{\left(  \sigma
-\delta\right)  \left(  \gamma+\sigma(1-q)\right)  }{\delta},$ $A_{0,0}%
^{(2)}\allowbreak=\allowbreak\frac{\gamma\delta}{\sigma\left(  \gamma
+\delta\right)  (\gamma+q\delta+(1-q)\sigma)}\left[  2\right]  _{q}%
\sqrt{\left(  \sigma-\delta\right)  \left(  \sigma+\gamma\right)  },$
$A_{0,1}^{(2)}\allowbreak=\allowbreak\frac{\gamma\delta}{\sigma\left(
\gamma+\delta\right)  (\gamma+q\delta+(1-q)\sigma)}\frac{\left(  \sigma
+\gamma\right)  \left(  q\delta+(1-q)\sigma\right)  }{\gamma},$ $A_{1,0}%
^{(2)}\allowbreak=-\left[  2\right]  _{q}\frac{\gamma\delta(\sigma-\delta
)}{\sigma\left(  \gamma+\delta\right)  (\gamma+q\delta+(1-q)\sigma)}$ we get
immediately (\ref{linear}) (\ref{dual1}) and (\ref{dual2}). To get
(\ref{warboth}) we perform the following calculations. $\operatorname{var}%
\left(  X_{\sigma}|\mathcal{F}_{\leq\sigma-\delta,\geq\sigma+\gamma}%
^{X}\right)  \allowbreak=\allowbreak\mathbb{E}\left(  X_{\sigma}%
^{2}|\mathcal{F}_{\leq\sigma-\delta,\geq\sigma+\gamma}^{X}\right)
\allowbreak-\allowbreak\left(  \mathbb{E}\left(  X_{\sigma}|\mathcal{F}%
_{\leq\sigma-\delta,\geq\sigma+\gamma}^{X}\right)  \right)  ^{2}$. Hence
\begin{align*}
\operatorname{var}\left(  X_{\sigma}|\mathcal{F}_{\leq\sigma-\delta
,.\geq\sigma+\gamma}^{X}\right)  \allowbreak &  =\allowbreak\frac{\delta
\gamma}{\left(  \sigma\left(  1-q\right)  +\gamma+q\delta\right)  }%
\allowbreak+\allowbreak X_{\sigma-\delta}^{2}(\frac{\gamma\left(  \left(
1-q\right)  \sigma+\gamma\right)  }{\left(  \delta+\gamma\right)  \left(
\sigma\left(  1-q\right)  +\gamma+q\delta\right)  }\allowbreak-\allowbreak
\frac{\gamma^{2}}{\left(  \delta+\gamma\right)  ^{2}})\allowbreak\\
&  +\allowbreak X_{\sigma+\gamma}^{2}(\frac{\delta\left(  \left(  1-q\right)
\sigma+q\delta\right)  }{\left(  \delta+\gamma\right)  \left(  \sigma\left(
1-q\right)  +\gamma+q\delta\right)  }\allowbreak-\allowbreak\left(
\frac{\delta}{\delta+\gamma}\right)  ^{2})\allowbreak\\
&  +\allowbreak X_{\sigma-\delta}X_{\sigma+\gamma}(\frac{\left(  q+1\right)
\delta\gamma}{\left(  \delta+\gamma\right)  \left(  \sigma\left(  1-q\right)
+\gamma+q\delta\right)  }\allowbreak-\allowbreak\frac{2\gamma\delta}{\left(
\delta+\gamma\right)  ^{2}}).
\end{align*}
After some simplifications we get
\begin{gather*}
\operatorname{var}\left(  X_{\sigma}|\mathcal{F}_{\leq\sigma-\delta,\geq
\sigma+\gamma}^{X}\right)  \allowbreak=\allowbreak\frac{\delta\gamma}{\left(
\sigma\left(  1-q\right)  +\gamma+q\delta\right)  }\allowbreak+\allowbreak
\frac{\left(  1-q\right)  \delta\gamma\left(  \sigma+\gamma\right)  }{\left(
\gamma+\delta\right)  ^{2}\left(  \sigma\left(  1-q\right)  +\gamma
+q\delta\right)  }X_{\sigma-\delta}^{2}\allowbreak\\
+\allowbreak\frac{\left(  1-q\right)  \delta\gamma\left(  \sigma
-\delta\right)  }{\left(  \gamma+\delta\right)  ^{2}\left(  \sigma\left(
1-q\right)  +\gamma+q\delta\right)  }X_{\sigma+\gamma}^{2}\allowbreak
-\allowbreak\frac{\left(  1-q\right)  \gamma\delta\left(  2\sigma
-\delta+\gamma\right)  }{\left(  \gamma+\delta\right)  ^{2}\left(
\sigma\left(  1-q\right)  +\gamma+q\delta\right)  }X_{\sigma-\delta}%
X_{\sigma+\gamma}\allowbreak\\
=\allowbreak\frac{\delta\gamma}{\left(  \sigma\left(  1-q\right)
+\gamma+q\delta\right)  \left(  \gamma+\delta\right)  ^{2}}\allowbreak\\
\times\allowbreak(\left(  \gamma+\delta\right)  ^{2}\allowbreak+\allowbreak
\left(  1-q\right)  \allowbreak\left(  X_{\sigma+\gamma}-X_{\sigma-\delta
}\right)  \allowbreak\left(  \left(  \sigma-\delta\right)  X_{\sigma+\gamma
}\allowbreak-\allowbreak\left(  \sigma+\gamma\right)  X_{\sigma-\delta
}\right)  ).
\end{gather*}

\end{proof}

\begin{proof}
[Proof of Remark \ref{kowariancje}]We have:
\begin{align*}
\operatorname{cov}\left(  \left(  X_{\tau}-X_{\sigma}\right)  ^{2},\left(
X_{\omega}-X_{\upsilon}\right)  ^{2}\right)  \allowbreak &  =\allowbreak
\newline\mathbb{E}\left(  \left(  X_{\tau}-X_{\sigma}\right)  ^{2}%
\mathbb{E(}\left(  X_{\omega}-X_{\upsilon}\right)  ^{2}|\mathcal{F}_{\leq\tau
}^{X})\right)  \allowbreak-\allowbreak(\omega-v)(\tau-\sigma)\allowbreak
=\allowbreak0,\\
\operatorname{cov}\left(  \left(  X_{\tau}-X_{\sigma}\right)  ^{3},\left(
X_{\omega}-X_{\upsilon}\right)  ^{3}\right)  \allowbreak &  =\allowbreak
-\left(  1-q\right)  \left(  \omega-\upsilon\right)  \mathbb{E}X_{\upsilon
}\left(  X_{\tau}-X_{\sigma}\right)  ^{3}\allowbreak\\
&  =-\left(  1-q\right)  \left(  \omega-\upsilon\right)  \mathbb{E}X_{\tau
}\left(  X_{\tau}-X_{\sigma}\right)  ^{3}\allowbreak\\
&  =\allowbreak-\left(  1-q\right)  \left(  \omega-\upsilon\right)
\mathbb{E}X_{\sigma}\left(  X_{\tau}-X_{\sigma}\right)  ^{3}\\
-\left(  1-q\right)  \left(  \omega-\upsilon\right)  \mathbb{E}\left(
X_{\tau}-X_{\sigma}\right)  ^{4}  &  =-\left(  1-q\right)  \left(  \tau
-\sigma\right)  \left(  \omega-\upsilon\right)  \left(  \tau\left(
2+q\right)  -\sigma\left(  1+2q\right)  \right)  .
\end{align*}
$\allowbreak$
\end{proof}

\begin{proof}
[Proof of Corollary \ref{rozk_q_wienera}]We multiply both sides of both
formulae given in assertion 4 of the Theorem \ref{rozk_q_wienera}) by
$s^{n}/\left[  n\right]  _{q}$ and sum over $n$ from $0$ to $\infty$. Now one
has to use that fact that such a sum is absolutely convergent in the case of
(\ref{prop1}) for $s^{2}\sigma(1-q)<1,$ while in the case of (\ref{prop2}) for
$s^{2}(1-q)/\sigma<1.$Changing order of summation and conditional expectation
and applying formula 1.4 of \cite{bms} giving generating function of $q-$
Hermite polynomials, leads to (\ref{general_martyngal}) and
(\ref{inverse_martyngal}) respectively. For $q=1$ we use the same formulae but
this time multiplied by $s^{n}/n!$ and after summing up with respect to $n$ we
apply well known formula for generating functions of Hermite polynomials.
\end{proof}

\end{document}